\documentclass[11pt]{article}

\usepackage{amssymb,amsmath,amsfonts,amssymb}
\usepackage{graphics,graphicx,color}
\usepackage{empheq}
\usepackage{slashed}

\usepackage{tikz,tkz-fct}
\usetikzlibrary{math}

\def\@abssec#1{\vspace{.05in}\footnotesize \parindent .2in
{\bf #1. }\ignorespaces}

\graphicspath{{/EPSF/}{Figures/}}
\DeclareGraphicsExtensions{.eps}

\usepackage[margin=1.0in]{geometry}

\newcommand{\vvvert}{|\!\!\:\|}

\newtheorem{theorem}{Theorem}[section]
\newtheorem{lemma}[theorem]{Lemma}

\newtheorem{remark}[theorem]{Remark}

\def \Rm {\mathbb R}
\def \Nm {\mathbb N}
\def \Cm {\mathbb C}
\def \Zm {\mathbb Z}

\newcommand{\eps}{\varepsilon}

\newcommand{\dsum}{\displaystyle\sum}
\newcommand{\dint}{\displaystyle\int}

\newcommand{\pdr}[2]{\dfrac{\partial{#1}}{\partial{#2}}}

\newcommand{\dr}[2]{\dfrac{d{#1}}{d{#2}}}

\newcommand{\aver}[1]{\langle {#1} \rangle}

\newcommand{\tx}{{\tilde x}} \newcommand{\ty}{{\tilde y}}
\newcommand{\tal}{{\tilde\alpha}}

\newcommand{\mE}{\mathfrak E}  

\newcommand{\mL}{\sqrt{\mathcal E}} 

\newcommand{\mS}{\mathcal S} \newcommand{\mT}{\mathcal T}
 \newcommand{\mX}{\mathcal X}
\newcommand{\rA}{{\rm A}} 
\newcommand{\rL}{{\rm L}}  \newcommand{\rM}{{\rm M}} \newcommand{\rT}{{\rm T}}
\newcommand{\rx}{{\rm x}} \newcommand{\ry}{{\rm y}}

\newcommand{\fa}{{\mathfrak a}}

\newcommand{\mCS}{\mathbb S}

\newcommand{\cout}[1]{}

\newcommand{\sgn}[1]{\,{\rm sign}(#1)}

\newcommand{\rJ}{{\rm J}}
\newcommand{\rP}{{\rm P}} \newcommand{\rPm}{{\rm P_{m}}}   \newcommand{\rPmz}{{\rm P_{m0}}}  
\newcommand{\rPmo}{{\rm P_{m1}}} \newcommand{\rPmt}{{\rm P_{m2}}}  
\newcommand{\trPm}{{\rm \tilde P_m}}   \newcommand{\rPmzi}{{\rm P_{m0}^{-1}}}  
\newcommand{\rPmj}{{\rm P_{mj}}} 
\newcommand{\rr}{{\rm r}}

\newcommand{\xim}{\xi_{\rm min}}
\newcommand{\epm}{\epsilon_m}

 \renewcommand{\arraystretch}{1.5}

\title{Semiclassical propagation along curved domain walls}

\author{Guillaume Bal \thanks{University of Chicago; {\tt guillaumebal@uchicago.edu}}}

\begin{document}
 
\maketitle


\begin{abstract}
  We analyze the propagation of two-dimensional dispersive and relativistic wavepackets localized in the vicinity of the zero level set $\Gamma$ of a slowly varying domain wall modeling the interface separating two insulating media. We propose a semiclassical oscillatory representation of the propagating wavepackets and provide an estimate of their accuracy in appropriate energy norms. We describe the propagation of relativistic modes along $\Gamma$ and analyze dispersive modes by a stationary phase method. In the absence of turning points, we show that arbitrary smooth localized initial conditions may be represented as a superposition of such wavepackets. In the presence of turning points, the results apply only for sufficiently high-frequency wavepackets. The theory finds applications both for Dirac systems of equations modeling topologically non-trivial systems as well as Klein-Gordon equations, which are topologically trivial.
 \end{abstract}

\renewcommand{\thefootnote}{\fnsymbol{footnote}}
\renewcommand{\thefootnote}{\arabic{footnote}}

\renewcommand{\arraystretch}{1.1}







%
\section{Introduction}

This paper concerns the propagation of two-dimensional localized wavepackets in a slowly varying domain wall that confines wave dynamics to the neighborhood of an interface. This interface may be interpreted as separating two insulating media. The slow variation of the underlying physical coefficients allows for a semiclassical description of the wavepacket propagation. Our objective is to describe how the geometry of the domain wall and  the other constitutive coefficients of the model influence wave propagation. 


Let $0<\eps$ be a small parameter and consider for concreteness the semiclassical Dirac equation
\begin{align}\label{eq:unpDirac}
  L_D = \eps D_t + \eps D_x \sigma_1 + \eps D_y \sigma_2 + \kappa(x,y) \sigma_3
\end{align}
with $D_a=-i\partial_a$ for $a\in\{t,x,y\}$, $t\in\Rm$ a time variable, $(x,y)\in\Rm^2$ two spatial variables, $\sigma_{1,2,3}$ the standard Pauli matrices and $\kappa(x,y)$ a mass term. Consider the case $\kappa(x,y)=y$ of a straight domain wall confining wavepackets in the neighborhood of $y=0$.  We observe by inspection that 
\begin{align}\label{eq:u0Dirac}
   u_0(t,x,y) = \frac{1}{\sqrt\eps} {\mathfrak f}\Big(\frac{x+t}{\sqrt\eps}\Big) \pi^{-\frac14} e^{-\frac12 (\frac y{\sqrt\eps}) ^2} \frac{1}{\sqrt2} \begin{pmatrix} 1  \\ -1 \end{pmatrix}
\end{align}
is a normalized solution of $L_Du=0$ in $L^2(\Rm^2;\Cm^2)$ for all $t\in\Rm$ when ${\mathfrak f}$ is normalized in $L^2(\Rm;\Cm)$. This expression describes a wavepacket in the $\sqrt\eps-$vicinity of a moving center $(-t,0)$ propagating with speed $-1$ along the $x-$axis. The $\sqrt\eps-$scaling is dictated by the interplay between the transport term $\eps D_y$ and the confining one $y\sigma_3$. 


As a spectral decomposition of $L_D$ demonstrates (see, e.g., \cite{B-EdgeStates-2018,bal2021edge} and the calculations in section \ref{sec:Dirac}), no similar wavepacket propagates with speed $+1$ along the $x-$axis. This reflects the non-trivial topology of the Dirac operator and the transport asymmetry along the interface $y=0$. This asymmetric transport along interfaces separating insulating materials is an important motivation for the analysis of \eqref{eq:unpDirac}.  Dirac equations are indeed ubiquitous macroscopic models in topological phases of matter 
\cite{moessner2021topological,WI}, 
in particular one-particle models of topological insulators \cite{volovik2009universe,Be13,PSB16,lee2019elliptic,B-BulkInterface-2018,B-EdgeStates-2018,bal2021edge,Drouot:21,drouot2020edge}; 
see also 
\cite{TKNN,avron1994,BES94,dombrowski2011quantization,quinn2022asymmetric} 
for integer quantum hall phases in the presence of a magnetic field not considered here.  A general principle called a bulk-edge correspondence states that transport along that interface has to be asymmetric with quantized asymmetry; see e.g. 
\cite{EG02,GP,PSB16,Drouot:19b,bal2022topological,bal2023topological,Drouot2020microlocal,bal2022multiscale,QB-NUMTI-2021} 
for mathematical analyses of the bulk-edge correspondence in different settings. 


Wavepackets such as \eqref{eq:u0Dirac} encode this asymmetry \cite[\S1.4]{bal2021edge} for flat interfaces. They persist in the presence of slowly varying (i.e., $\eps\ll1$) curved domain walls
\cite{bal2022magnetic,bal2021edge,hu2022traveling}, a result that we will re-establish below. However, for a domain wall $\kappa(x,y)=y+V(x,y)$ with $V$ a (compactly supported) perturbation and when $\eps=1$ (i.e., for a non-slowly varying wall), this asymmetric transport is encoded by a perturbation-dependent dispersive mode rather than the relativistic (i.e., non-dispersive) mode \eqref{eq:u0Dirac}; see, e.g., \cite{B-EdgeStates-2018,bal2023asymmetric}. While the bulk-edge correspondence implies that the transport asymmetry persists in the presence of perturbations, only in the semiclassical, weakly scattering, regime $\eps\ll1$, is that asymmetry described by a mode such as \eqref{eq:u0Dirac}.


Furthermore, the same spectral decomposition of $L_D$ shows that the domain wall $\kappa(x,y)=y$ confines many other modes in the vicinity of $y=0$. For instance, we observe that 
\begin{align}\label{eq:upm}
u_{1,\pm}(t,x,y)=\int_{\Rm} f_\pm(\xi) \frac{1}{\sqrt\eps} e^{i\frac{\mp\sqrt{2+\xi^2} t + \xi x}{\sqrt\eps}}  \frac{e^{-\frac12 (\frac y{\sqrt\eps}) ^2}}{\pi^{\frac14}}  (\sigma_1+\sigma_3) \begin{pmatrix} 1  \\  (\pm \sqrt{2+\xi^2}-\xi) \frac{y}{\sqrt\eps}\end{pmatrix} d\xi
\end{align} 
for $f_\pm(\xi)$ decaying sufficiently rapidly at infinity also belongs to the kernel of $L_D$.  This is a superposition (in the variable $\xi$) of plane waves with energy $E_\pm=\pm\sqrt{2+\xi^2}$. Such plane waves propagate with different group velocities, both positive and negative, $\partial_\xi E_\pm = \xi/E_\pm$ for different values of $\xi$ and therefore experience dispersion. 

Dispersive effects for such modes over long times are analyzed in detail in section \ref{sec:wae} by stationary phase methods \cite{dimassi1999spectral}. Denoting by $G_{1,\pm}(t,x,\xi)=\mp\sqrt{2+\xi^2} t + \xi x$ with $G_{1,\pm}(0,x,\xi)=x\xi$, we observe that $u_{1,\pm}(0,x,y)$ models a wavepacket concentrated in the $\sqrt\eps-$vicinity of $(x_0,y_0)=(0,0)$. As an application of stationary phase, we observe that for $x\in (-t,t)$, assuming $t>0$ to simplify, and defining $\xi_m$ as the unique solution to $\partial_\xi(\mp\sqrt{2+\xi^2} t + \xi x)=0$,  (i.e., $\xi_m=\pm(x/t) \sqrt2 (1-(x/t)^2)^{-\frac12}$) then $u_{1,\pm}(t,x,y)$ is well approximated by 
\begin{align*}
  f_\pm(\xi_m) \frac{1}{\eps^{\frac14}} e^{i\frac{G_{1,\pm}(t,x,\xi_m)}{\sqrt\eps}} e^{-i\frac \pi 4 \sgn{tG_\pm''}}  \Big( \frac{2\pi}{|tG_\pm''|}\Big)^{\frac12} \frac{e^{-\frac12 (\frac y{\sqrt\eps}) ^2}}{\pi^{\frac14}}  (\sigma_1+\sigma_3) \begin{pmatrix} 1  \\  (\pm \sqrt{2+\xi_m^2}-\xi_m) \frac{y}{\sqrt\eps}\end{pmatrix},
\end{align*}
where we have defined $G_\pm''=\partial^2_\xi G_{1,\pm}(t,x,\xi_m)$. Outside of the interval $x\in(-t,t)$, the solution $u_{1,\pm}(t,x,y)$ is asymptotically negligible compared to the above term as an application of non-stationary phase. 

In contrast to $u_0(t,x,y)$ in \eqref{eq:u0Dirac}, which is concentrated in the $\sqrt\eps$-vicinity of a {\em point} $(-t,0)$, the dispersive modes $u_{1,\pm}(t,x,y)$ are concentrated in the $\sqrt\eps$-vicinity of the {\em interval} $\{(x,0);\ -t\leq x\leq t \}$. Dispersion further causes the maximal amplitude of the wavepacket to be of order $\eps^{-\frac14} t^{-\frac12}$.

\medskip

Wave propagation with domain wall confinement is interesting in other contexts than the above Dirac system. Decompose $L_D=\eps D_t +H$ and formally square $\eps D_t=-H$ to obtain $\eps^2D_t^2=H^2$. Keeping one of the diagonal terms in this equation gives the following Klein-Gordon model
\begin{align}\label{eq:unpKG}
  L_{KG} =   \eps^2\partial_t^2 - \eps^2 \Delta  + y^2 - \eps = \eps^2\partial_t^2-\eps^2\partial_x^2 + \fa^*\fa
\end{align}
with $\fa=\eps\partial_y+y$. We then observe that the scalar functions
\begin{align}\label{eq:u0KG}
   u_\pm(t,x,y) = \frac{1}{\sqrt\eps} {\mathfrak f}\Big(\frac{x \pm t}{\sqrt\eps}\Big) \pi^{-\frac14} e^{-\frac12 (\frac y{\sqrt\eps}) ^2} 
\end{align}
are in the kernel of $L_{KG}$. This provides two relativistic modes confined by the domain wall $y^2$, or more generally by a domain wall $\kappa^2(x,y)$. Dispersive modes similar to \eqref{eq:upm} may also be constructed in the kernel of $L_{KG}$. 

The above wave equation with domain wall finds applications in the modeling of the Nambu-Goldstone modes that appear as a (static) domain wall partially breaks translational invariance, and propagate in a semiclassical picture along the zero-level set of $\kappa$, either a straight line when $\kappa(x,y)=y$, or the bottom of a `Mexican hat' when $\kappa(x,y)=x^2+y^2-1$, say; see \cite{watanabe2020counting} and its references for details on such ubiquitous modes. 

\medskip

The main objective of this paper is to describe the propagation of similar dispersive and relativistic (non-dispersive) wavepackets for general smooth domain walls $\kappa(x,y)$ and differential operators $L$ generalizing $L_D$ and $L_{KG}$. We assume the existence of a domain wall $\kappa(x,y)$ confining wavepackets in the vicinity of its zero level set $\Gamma=\kappa^{-1}(0)$, assumed to be a smooth one-dimensional non-intersecting curve in $\Rm^2$. 
Locally, generalizations of the above plane waves are approximately in the kernel of $L$ and modulations of such plane waves allow us to approximate solutions of $Lu=0$ with arbitrary accuracy in powers of $\eps$ for large (but not infinite) times in many favorable cases.

\medskip

With the Dirac and Klein-Gordon models serving as illustrations, we develop an abstract framework to analyze wave packet propagation along smooth slowly varying curved domain walls. This is detailed in section  \ref{sec:wae}. The verification of all required hypotheses is carried out for general Dirac operators in section \ref{sec:Dirac} and Klein-Gordon operators in section \ref{sec:KG}.  The main results and principal steps of the derivation presented in section \ref{sec:wae} are summarized as follows.

\medskip\noindent
1) {\em Rectification}. In the vicinity of $(x_0,y_0)\in\Gamma=\kappa^{-1}(0)$, the domain wall is well approximated by $(x-x_0,y-y_0)\cdot\nabla\kappa(x_0,y_0)$. This forms a linear domain wall similar to the one introduced in \eqref{eq:unpDirac} when $\kappa(x,y)=y$. It is convenient to work in a system of coordinates where the curve  $\Gamma$ is parametrized by $\tx$ and the linearization of the domain wall is of the form $\ty\partial_\ty\kappa$.  In such a system, the operator $\fa=\eps\partial_y+y$ naturally appears and explains why Hermite functions are ubiquitous in \eqref{eq:upm} and \eqref{eq:u0KG}. This change of variables, described in section \ref{sec:rectification}, is implemented by a map $(x,y)=\Phi(\tx,\ty)$ parametrizing a tubular neighborhood of $\Gamma$. This map is a diffeomorphism from a vicinity of $\ty=0$ to its image when $\Gamma$ is open (and then unbounded) and acts as a covering map when $\Gamma$ is closed and then diffeomorphic to a circle; see Figure \ref{fig:geom}.

The differential operator $L$ of interest written in the rectified coordinates is $P=\Phi^* L \Phi^{-*}$; see section \ref{sec:rectification} for notation. Solutions of $Lu=0$ are then analyzed in a $2\eta-$ neighborhood of $\{\ty=0\}=\Phi^{-1}(\Gamma)$ for $\eta>0$ chosen sufficiently small (independently of $\eps$) by means of $Pv=0$ with $v=\Phi^*u$ a pulled back version of $u$ in the rectified coordinates.
 
 \medskip\noindent
 2) {\em Wavepacket multiscale (WKB) ansatz}. Operators such as $D_x\sigma_1+D_y\sigma_2$ or $-\Delta$ on $\Rm^2$ admit (absolutely) continuous spectrum parametrized for instance by a two-dimensional Fourier variable. In the presence of a domain wall $y$ or $y^2$, the spectrum heuristically becomes discrete in the variable $y$ while it remains continuous in the orthogonal variable $x$. This leads us to the analysis of a countable number of branches of simple absolutely continuous spectrum. The latter are then parametrized by the following ansatz:
 \begin{align}\label{eq:vansatz}
   v(t,x,y) = \frac{1}{\sqrt\eps} \dint_{\Xi} e^{\frac{i}{\sqrt\eps} G_m(t,x,\xi;x_0)} \psi(t,x,\xi,\frac{y}{\sqrt\eps};x_0) d\xi,
\end{align}
where $m$ labels the branch while $\xi\in\Xi\subset\Rm$ parametrizes it. The function $\psi(t,x,\xi,z;x_0)$ involves `macroscopic' variables $(t,x)$ and `microscopic' variables $(\xi,z)$. The parameter $x_0$ is a reminder that wavepackets are localized in the vicinity of $(x_0,0)$ at time $t=0$. Finally, $G_m(t,x,\xi;x_0)$ is a phase term modeling the wave dispersion relations. It takes the form $G_{(0,-)}=\xi(t+x)$ or $G_{(0,\pm)}=\xi(\mp t+x)$ for the above relativistic modes and $G_{(1,\pm)}= \mp\sqrt{2+\xi^2} t + \xi x$ for the dispersive modes in \eqref{eq:upm}.

\medskip\noindent
 3) {\em Multiscale asymptotic expansion and wavepacket construction}. Looking for approximations of $Pv=0$ is performed by looking at $m$ fixed for $\rPm\psi=0$ with a multiscale operator $\eps^{-\frac q2}\rPm=\rPmz+\sqrt\eps \rPmo+\eps \rPmt$ (with $q=1$ for Dirac and $q=2$ for Klein-Gordon). Decomposing $\psi=\sum_{j\geq0}\eps^{\frac j2}\psi_j$ provides a sequence of formal equations that are solved in turn. As is standard in two-scale analyses, the first two equations characterize the leading approximation to $\psi$. 
 
The first equation $\rPmz\psi_0=0$, describes the local, microscopic equilibrium. It provides a Hamilton-Jacobi (eikonal) equation for $G_m$ appearing in \eqref{eq:vansatz} and the structure of the domain wall confinement. The second equation $\rPmz\psi_1+\rPmo\psi_0=0$ describes transport along the interface. This is analogous to standard semiclassical analyses in the absence of a domain wall \cite{dimassi1999spectral,zworski2012semiclassical} with the main difference that transport is now mainly confined along a one-dimensional interface.

Hamilton-Jacobi equations are typically solvable only for short times \cite{dimassi1999spectral,zworski2012semiclassical}; see \cite{drouot2022semiclassical} for a detailed analysis in the presence of similar confinement to the one we consider here. Since we aim at a construction of wavepackets that applies to long times, we look for solutions of the form
\begin{align}\label{eq:ansatzGm}
   G_m(t,x,\xi;x_0) = - B_m(\xi;x_0) t + A_m(x,\xi;x_0)
\end{align}
where $A_m$ and $B_m$ are appropriate functions with $\partial_x A_m(x_0,\xi;x_0)=\xi$ as expected so that the ansatz \eqref{eq:vansatz} approximately corresponds to a (semiclassical) Fourier transform at time $t=0$. We will obtain such solutions to the Hamilton Jacobi equations for a large class of wavepackets (large choice of the wavenumber support $\Xi$ in \eqref{eq:vansatz}), and in fact for all possible reasonable wavepackets (i.e., with $\Xi=\Rm$ or $\Xi=\{|\xi|\geq \xi_0\}$) for Dirac and Klein-Gordon operators when some appropriate coefficients are constant. Addressing general wavepackets when such coefficients are not constant would require the construction of solutions $G_m$ with initial conditions given by $G_m(0,x,\xi;x_0)=(x-x_0)\xi$  typically valid only for short times  \cite{drouot2022semiclassical,dimassi1999spectral,zworski2012semiclassical} as caustics may develop. Long-time or tunneling effects in the presence of such caustics are then complicated even in one-dimensional settings \cite{lindblad2020modified,ramond1996semiclassical} and we do not consider them here.
 
\medskip\noindent
4) {\em Main results on wavepacket analysis and dispersive effects}.  Once explicit approximate wavepackets $\psi^J_m=\sum_{j=0}^J \eps^{\frac j2}\psi_j$ for an arbitrary $J$ are constructed in the rectified coordinates, we may construct $v^J_m$ using \eqref{eq:vansatz} and then $u^J:=u^J_m=\Phi^{-*} v^J_m$ (see \eqref{eq:ufromv} for a more precise statement when $\Gamma$ is closed). We then obtain the following results:

(1) The first main result of the paper is Theorem \ref{thm:localerror} below, stating that the solution $u$ of $Lu=0$ with the same initial conditions as $u^J$ is such that $u-u^J$ is small in an appropriate energy norm. This requires the multi-scale expansion of $L$ to satisfy conditions (i)-(v) in section \ref{sec:mo} and the energy functional to satisfy condition (vi) in section \ref{sec:wpa}.

(2) The above initial condition for $u$ concentrates on the $m$th branch of continuous spectrum. Assuming that \eqref{eq:ansatzGm} holds for a complete set of branches $m$  and all $\xi\in\Xi=\Rm$ with real-valued phases $G_m$ (see assumptions (vii)-(viii) in section \ref{sec:wpa}), then we show that any sufficiently smooth initial condition for $u$ concentrated in the vicinity of $(x_0,y_0)\in\Gamma$ can indeed be represented by a superposition over the branches of continuous spectrum of explicit wavepackets of the form $\Phi^{-*} v^J$ with $v^J$ as in \eqref{eq:vansatz}.  These results, which essentially construct a semiclassical parametrix for $L$, are described in Theorem \ref{thm:globalerror}.

(3) Once wavepackets of the form $\Phi^{-*} v$ have been constructed for a fixed branch $m$ with $v$ represented by \eqref{eq:vansatz}, we can analyze the leading term in the expansion in powers of $\eps$. Theorem \ref{thm:nondispersive} then generalizes the description of relativistic solutions such as \eqref{eq:u0Dirac} to curved domain walls.  The explicit expression of this mode for Dirac equations is given in Theorem \ref{thm:propzero}. This extends results obtained in \cite{bal2021edge} for the Dirac equation using a different ansatz. 

(4) The leading term of dispersive modes requires a more detailed analysis. Under assumptions of non-vanishing curvature (to promote dispersion) of $G_m$ stated in condition (ix) in section \ref{sec:wpa}, we analyze the dispersive effects of superpositions of generalized plane waves of the form \eqref{eq:upm} by means of stationary phase in Theorem \ref{thm:dispersive}. 

All approximation results in the above theorems are valid up to times  $\rT$ (essentially) of order $\eps^{\frac12} \rT\ll1$. At larger times, the curved domain wall allows for the coupling of modes belonging to different branches $m$ by scattering. An ansatz of the form \eqref{eq:vansatz} ceases to be valid and we leave the semiclassical regime.

 
\medskip
 
While the results of  section \ref{sec:wae} are applied to Dirac equations in section \ref{sec:Dirac} and Klein-Gordon operators in section \ref{sec:KG}, we expect them to be valid for a broader class of models of wave propagation such as those, e.g., in \cite{volovik2009universe} with applications in topological superconductors,  in \cite{bal2022multiscale} with applications to Floquet topological insulators, as well as the asymmetric propagation of shallow water waves along the equator confined by the Coriolis force \cite{delplace2017topological,souslov2019topological,bal2023topological}.
 
We also restrict ourselves to a two-dimensional setting with one-dimensional interface. There should be no fundamental obstruction, beyond technical difficulties, to looking at propagation along one-dimensional curves obtained as the intersection of $d-1$ domain walls  in $\Rm^d$ (with $v(t,x,y)$ still given by \eqref{eq:vansatz} now for $y\in\Rm^{d-1}$) as they appear in the topological classification and analysis of high-order topological insulators \cite{bal2023topological}.
\section{Wavepacket construction and analysis}\label{sec:wae}
\subsection{Domain wall and rectified problem}\label{sec:rectification}

This section introduces a change of variables $\Phi$ rectifying the level set $\Gamma=\kappa^{-1}(0)$. 
The domain wall is modeled by a smooth scalar, bounded, function $\kappa(x,y)\in C^\infty(\Rm^2)$ whose zero-level set $\Gamma=\kappa^{-1}(0)$ is assumed to be connected. We impose the condition 
\begin{align}\label{eq:condkappa}
\inf_{(x,y)\in\Gamma} |\nabla\kappa(x,y)|>0
\end{align}
to avoid any degeneracy along $\Gamma$. With this assumption, $\Gamma$ is either an open infinite curve or a closed compact, non-intersecting, curve diffeomorphic to the unit circle. The construction and analysis of semiclassical wavepackets propagating in the vicinity of $\Gamma$ is significantly simplified when moving to rectified coordinates, where $\Gamma$ becomes the $x-$axis when it is open and a sub-interval of the $x-$axis when it is closed.

The zero level set $\Gamma=\kappa^{-1}(0)$ is parametrized by arclength along $\Gamma$: $\tx\to \Phi(\tx,0)$ such that the unit tangent vector to $\Gamma$ is given by $\tau(\tx)=\partial_\tx \Phi(\tx,0)$ with $|\tau(\tx)|=1$. We denote by $\nu(\tx)$ the normal vector so that $(\tau(\tx),\nu(\tx))$ is an orthonormal basis of $\Rm^2$ with positive orientation.

\begin{figure}[ht!] 
\begin{tikzpicture}[scale=2]
\begin{scope}[shift={(0.2,0)}]
\clip (3.5,-.5)  rectangle (9,2)  plot [smooth cycle, tension=0.6] coordinates {(4.4,0.4)   (5,0.2)  (5.8,0.6) (6.5773,0.5421)(6.4905,1.1074)  (5.9752,1.2828) (5.4,1.4) (4.6,1) };
\draw[line width=1cm] plot [smooth cycle, tension=0.6] coordinates {(4.4,0.4) (5,0.2) (5.8,0.6) (6.5773,0.5421)(6.4905,1.1074)  (5.9752,1.2828) (5.4,1.4) (4.6,1) };
\draw[line width=.975cm,white] plot [smooth cycle, tension=0.6] coordinates {(4.4,0.4) (5,0.2) (5.8,0.6) (6.5773,0.5421)(6.4905,1.1074)  (5.9752,1.2828) (5.4,1.4) (4.6,1) };
\draw[line width=.52cm] plot [smooth cycle, tension=0.6] coordinates {(4.4,0.4) (5,0.2) (5.8,0.6) (6.5773,0.5421)(6.4905,1.1074)  (5.9752,1.2828) (5.4,1.4) (4.6,1) };
\draw[line width=.42cm,white] plot [smooth cycle, tension=0.6] coordinates {(4.4,0.4) (5,0.2) (5.8,0.6) (6.5773,0.5421)(6.4905,1.1074)  (5.9752,1.2828) (5.4,1.4) (4.6,1) };
\draw[line width=.125mm] plot [smooth cycle, tension=0.6] coordinates {(4.4,0.4) (5,0.2) (5.8,0.6) (6.5773,0.5421)(6.4905,1.1074)  (5.9752,1.2828) (5.4,1.4) (4.6,1) };
\draw[line width=1cm,dash pattern=on 0.5pt off 20pt] plot [smooth cycle, tension=0.6] coordinates {(4.4,0.4) (5,0.2) (5.8,0.6) (6.5773,0.5421)(6.4905,1.1074)  (5.9752,1.2828) (5.4,1.4) (4.6,1) };
\end{scope}

\node at (3.4,1.4) {{\huge $\Phi$}};
\draw [line width = .8mm,->] (3,1.1) -- (4,1.1);
\node at (6.1,1.9) {{\Large $\Gamma=\{\kappa=0\}$}};
\node at (6.2,1.0) {{\Large $\kappa<0$}}; \node at (0,.2) {{\Large $\kappa<0$}};
\node at (4.65,1.55)  {{\Large $\kappa>0$}};\node at (0,1.5)  {{\Large $\kappa>0$}};
\node at (1.4,.25) {{\Large $(\tilde x,\tilde y)$}};\node at (6.2,-.2) {{\Large $(x,y)=\Phi(\tilde x,\tilde y)$}};
\node at (1.365,.5) {$\bullet$};\node at (5.25,-.05) {$\bullet$};
\node at (4.465,.5) {$\bullet$}; \draw [line width = .4mm,->] (4.465,.5) -- (4.465,1); \draw [line width = .4mm,->] (4.465,.5) -- (3.965,.5)  ;
\node at (4.4,1.1) {{\large $\tau$}};  \node at (4.05,.34) {{\large $\nu$}};

\begin{scope}[shift={(-3,.5)}, scale=0.660]
\draw[line width=.8mm] (4,0.4) -- (8,0.4);
\draw[line width=.4mm] (4,0) -- (8,0);
\draw[line width=.4mm] (4,.8) -- (8,0.8);
\draw[line width=.4mm] (4.2,0)-- (4.2,.8);\draw[line width=.4mm] (5,0)-- (5,.8);\draw[line width=.4mm] (5.8,0)-- (5.8,.8);
\draw[line width=.4mm] (6.6,0)-- (6.6,.8);\draw[line width=.4mm] (7.4,0)-- (7.4,.8);
\draw[dashed] (3.5,0) -- (8.5,0);
\draw[dashed] (3.5,0.4) -- (8.5,0.4);
\draw[dashed] (3.5,0.8) -- (8.5,0.8);
\end{scope}
\end{tikzpicture}
\label{fig:geom}\caption{Geometry of the domain wall $\kappa$ with tubular neighborhood of a closed $\Gamma=\kappa^{-1}(0)$ and rectification mapping $\Phi$.}
\end{figure}
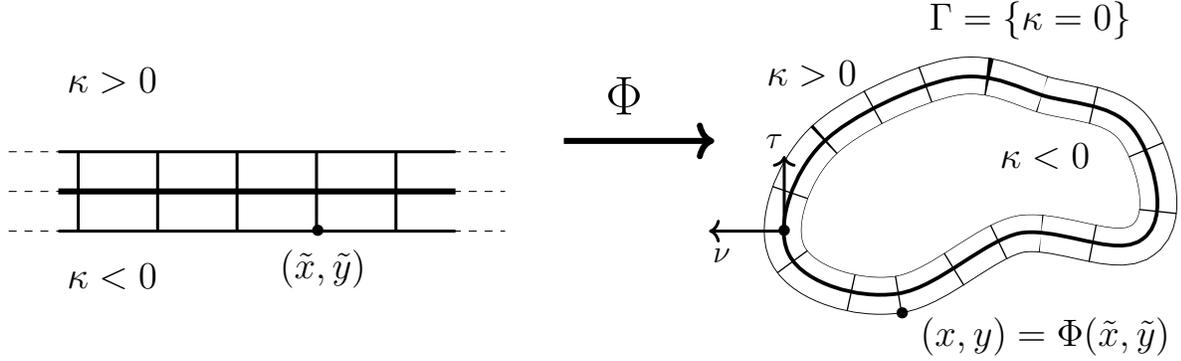



Let $\eta>0$. We then define on $\Rm\times [-2\eta,2\eta]$ the map
\begin{align}\label{eq:rectification}
   (\tx,\ty) \to (x,y) = \Phi(\tx,\ty) := \Phi(\tx,0) + \ty \nu (\tx).
\end{align} 
For  $\eta$ sufficiently small, by assumption on $\kappa$, $\Phi$ maps $\Rm\times [-2\eta,2\eta]$ to a tubular neighborhood of $\Gamma$. The variables $(\tx,\ty)$ will be referred to as the rectified coordinates. The vector $\tau$ above is defined up to a choice of orientation. We choose $(\tau,\nu)$ such that $\kappa\circ \Phi(\tx,\ty) >0$ when $\ty>0$. 

We state several standard results on the geometry of curves in a plane.
When $\Gamma$ is open, then $\Phi$ is a diffeomorphism from $\Rm\times[-2\eta,2\eta]$ to its image. When $\Gamma$ is closed, and hence compact, it has Euclidean length $\rL<\infty$.  Then $\Phi$ is a covering map (with fiber $\Zm$) in the sense that all points $(\tx+k\rL,\ty)$ for $k\in\Zm$ are mapped to the same point $\Phi(\tx,\ty)$ near $\Gamma$.

The Jacobian of the transformation  is defined as $\rJ(\tx,\ty)={\rm D}\Phi(\tx,\ty)$. We verify that it is given by $((1+\varkappa(\tx)\ty)\tau(\tx) \ \nu(\tx))$ as a matrix of column vectors, where $\varkappa(\tx)$ is the curvature defined as the scalar function such that $\partial_\tx \nu (\tx) = \varkappa(\tx) \tau(\tx)$. Thus, 
\begin{align}\label{eq:detJ}
    \det\,\rJ = 1+\varkappa(\tx)\ty>0 \quad \mbox{ on } \quad \Rm\times[-2\eta,2\eta],
\end{align}
provided $\eta$ is chosen small enough, which we assume from now on.

\medskip

Associated to the mapping $\Phi$ is a linear transformation of functions, the pull-back $\Phi^*$ defined for each smooth scalar function $f(x,y)$ by $\Phi^*f(\tx,\ty)=f\circ\Phi(\tx,\ty)$. The inverse pull-back $\Phi^{-*}=(\Phi^{-1})^*$ is defined globally when $\Gamma$ is open since then $\Phi$ is a diffeomorphism and is defined locally (for functions with sufficiently small support) when $\Gamma$ is closed. When $\Gamma$ is closed, then $\Phi^{-*}$ is defined globally on periodic functions. For a function $v(\tx,\ty)$ with sufficiently rapid decay as $|\tx|\to\infty$, we define the periodization
\begin{align}\label{eq:periodization}
  v_\sharp(\tx,\ty) = \dsum_{k\in\Zm} v(\tx+k\rL,\ty).
\end{align}
On such periodic functions $v_\sharp(\tx,\ty)$, the inverse pull-back $u=\Phi^{-*}v_\sharp = v_\sharp \circ\Phi^{-1}$ is globally defined.

\medskip

The above pull-back operator is convenient to pull back a differential operator in the $(x,y)$ coordinates to one in the rectified $(\tx,\ty)$ coordinates. In particular, since differentiation is local, we define the pulled-back derivatives
\begin{align}\label{eq:pullbackD}
   \tilde D = \Phi^* D \Phi^{-*} = \rJ^{-T} D
\end{align}
where we use the notation $A^{-T}=(A^{-1})^T$ with $T$ transposition. More generally, for a differential operator $L$, we define the pulled-back differential operator 
\begin{align}\label{eq:pullbackL}
  P=\Phi^* L \Phi^{-*}.
\end{align}
An operator $Lu=0$ for $u$ supported in the image of $\Phi$ will therefore locally be equivalent to a problem $Pv=0$ where $v=\Phi^* u$, or $u=\Phi^{-*}v$.

\begin{remark}[Notation]\label{rem:notationscales}\rm
To simplify expressions, the rectified variables $(\tx,\ty)$ are also denoted by $(x,y)$ when no confusion is possible. 
We systematically use $u(t,x,y)$ for a function in the original, unrectified coordinates, $v(t,x,y)$ (instead of $v(t,\tilde x,\tilde y)$) for a function in the rectified coordinates with $v=\Phi^* u$, $\Psi(t,x,z)$ a function in the rescaled rectified coordinates with $z=\frac y{\sqrt\eps}$, and finally $\psi(t,x,\xi,z)$ as a  function in the multiscale rescaled rectified coordinates (as in \eqref{eq:vansatz}).
\end{remark}

\subsection{Wavepacket ansatz} \label{sec:wa}

We wish to construct explicit approximations of solutions to the semiclassical differential equation
\begin{align}\label{eq:Lu}
  L(x,y,\eps D_t,\eps D_x,\eps D_y) u(t,x,y) =0.
\end{align}
Here, $(x,y)\in\Rm^2$, $D_a=-i\partial_a$ for $a\in \{t,x,y\}$, and $L(x,y,p_t,p_x,p_y)$ is a polynomial in the last three components with smooth $(x,y)-$dependence acting on functions $u(t,x,y)$ with values in $\Cm^l$, where $l=2$ for Dirac and $l=1$ for Klein-Gordon models. 
$0<\eps\leq1$ is a semiclassical parameter. We aim to construct solutions of $Lu=0$ asymptotically as $\eps\to0$. 

Using \eqref{eq:pullbackD}, we introduce the  pulled back operator  $P=\Phi^{*} L \Phi^{-*}$ such that for $v=\Phi^{*}u$, \eqref{eq:Lu} is equivalent in a neighborhood of $\Gamma$ to
\begin{align}\label{eq:PPsi}
  P(x,y,\eps D_t,\eps D_x,\eps D_y) v(t,x,y) =0.
\end{align}
As mentioned in Remark \ref{rem:notationscales}, we drop the $\ \tilde{}\ $ in the rectified coordinates.

For $\eta>0$, denote by $0\leq \chi_\eta(y)\leq 1$ a $C_c^\infty(\Rm)$ function equal to $1$ for $|y|\leq\eta$ and equal to $0$ for $|y|\geq 2\eta$.   We then define the operator
\begin{align*}
   P_\eta = P \chi_\eta = \Phi^{*} L \Phi^{-*} \chi_\eta.
\end{align*}
With a slight abuse of notation, the operator $P_\eta$ thus applies to functions defined on $(t,x,y)\in [-\rT,\rT]\times\Rm^2$ rather than only in a vicinity of $\Phi^{-1}(\Gamma)$.

The wavepacket construction hinges on a scale separation in the rectified $(x,y)$ variables. We introduce the rescaled transverse variable $z$ such that $y=\sqrt\eps z$. While the domain wall prevents propagation away from $\Gamma$, propagation along $\Gamma$ is allowed and modeled by branches of continuous spectrum. This is the motivation for the following ansatz
\begin{align}\label{eq:Psi}
   v(t,x,y) & =  \frac{1}{\sqrt\eps} \dint_{\Xi} e^{\frac i{\sqrt\eps} G_m(t,x,\xi)} \psi \big(t,x,\xi,\frac{y}{\sqrt\eps}\big) d\xi\\ \nonumber
   \eps^{\frac14} v(t,x,\sqrt\eps z) =: \Psi(t,x,z) &= \frac{1}{\eps^{\frac14}} \dint_{\Xi} e^{\frac i{\sqrt\eps} G_m(t,x,\xi)} \psi (t,x,\xi,z) d\xi.
\end{align}
The above scalings are introduced so that all functions are normalized to $O(1)$ in the $L^2$ sense in their respective spatial variables. Here, $\Xi$ is an open subset of $\Rm$ that in the applications will take the form $\Xi=\Rm$ or $\Xi=\{\xi\in\Rm; \ |\xi|> \xim\}$. The role of $\xim$ is to ensure that wavepackets with support in $\Xi$ do not encounter caustics (turning points).

Above, $G_m(t,x,\xi;x_0)$ is a phase function depending on a  variable $\xi$ parametrizing the $m$th branch of continuous spectrum, where the index $m$ takes a countable number of values. The phases $G_m$ depend on a parameter $x_0$  capturing the location of the initial wavepacket, although the explicit dependence in $x_0$ will often be omitted. The phase functions we consider are real-valued and will be homogeneous of degree $1$ in the variables $(\xi,\sqrt m)$ properly defined.  The wavefront $e^{\frac i{\sqrt\eps} G_m(t,x,\xi)} \psi \big(t,x,\xi,\frac{y}{\sqrt\eps}\big)$ is then similar to a standard WKB ansatz \cite{dimassi1999spectral,zworski2012semiclassical} with the property that it describes propagation along the interface (variables $(t,x,\xi)$) and confinement across it (variables $(m,y)$). The superposition in $\xi$ generates spatially localized wavepackets $v(t,x,y)$ whose evolution we aim to describe.

Denote by $\tilde P_\eta=P_\eta(x,\sqrt\eps z,\eps D_t,\eps D_x,\sqrt\eps D_z)$ the (pulled-back) operator (of $P_\eta$) in the rescaled variable $z\to y(z) =\sqrt\eps z$.
Then we have
\begin{align}\nonumber
  \eps^{\frac14} P_\eta v(t,x,\sqrt\eps z) = \tilde P_\eta\Psi(t,x,z) = & \frac{1}{\eps^{\frac14}} \dint_{\Xi} e^{\frac i{\sqrt\eps} G_m(t,x,\xi)} \rP_m\psi (t,x,\xi,z) d\xi,
   \\ 
  \label{eq:rPm}
  \rP_m(x,z ,\eps D_t,\eps D_x,\sqrt\eps D_z) = &\  P_\eta\big(x,\sqrt \eps z, \eps D_t + \sqrt\eps \partial_t G_m,\eps D_x + \sqrt\eps \partial_x G_m, \sqrt\eps D_z\big). 
\end{align}
We do not retain the explicit dependence of $\rPm$ on $\eta$ to simplify.


%
\subsection{Multiscale operator expansion} \label{sec:mo}
In \eqref{eq:rPm}, $\rP_m$ is the {\em multiscale} operator associated to $P_\eta$ for {\em branch} $m$.  By assumption, $P_\eta$ is a polynomial of degree $q$ in the last three components. The dependence of $P_\eta$ on $\sqrt\eps z$ involves the smooth domain wall $\kappa(x,\sqrt\eps z)$, which vanishes at $z=0$ and is approximately equal to its Taylor expansion in the vicinity of $z=0$. Combining these two observations, we assume the expansion:
\begin{align}\label{eq:expP}
  \rPm (x,z ,\eps D_t,\eps D_x,\sqrt\eps D_z)= \eps^{\frac q2} \big( \rPmz+\sqrt\eps \rPmo +\eps \rPmt \big).
\end{align}
The operators $\rPmz$ and $\rPmo$ involve $\eps-$independent coefficients and the leading term satisfies
\begin{align}\label{eq:P0}
  \rPmz = P_0 (x,z,\partial_t G_m, \partial_x G_m,D_z)
\end{align}
for a function $P_0$ that we {\em assume} is a polynomial of degree $q$ in the {\em last four} variables (with matrix-valued, $x-$dependent coefficients) and hence defined for $z\in\Rm$. We will obtain such a decomposition for Dirac and Klein-Gordon operators in the following sections based on the existence of a domain wall. In what follows, what matters is the properties of the operators $\rPmj$ independently of the presence of a domain wall or not.

The leading-order equation 
\begin{align} \label{eq:eikonal}
  \rPmz \psi =    P_0(x,z,\partial_t G_m,\partial_x G_m,D_z) \psi =0
\end{align}
is an {\em eikonal} equation (for $G_m$ and $\psi$). For each fixed $(x_0,0)\in\Phi^{-1}(\Gamma)$ about which the initial wavepacket is localized, we {\em assume} that it admits a countable family of solutions parametrized by $m\to (G_m(x,\xi;x_0),\phi_m(x,\xi,z;x_0))$, which holds in the cases of interest to us. 

\medskip

For each fixed $m$, we aim to construct accurate solutions of the multiscale equation $\rPm\psi=0$. Using \eqref{eq:Psi}, this constructs approximations $v(t,x,y)$ of $Pv=0$ and then moving back to the physical variables, approximations $u=\Phi^{-*}v$ of the original problem $Lu=0$.
Formally decomposing $\psi=\sum_{j\geq0} \eps^{\frac j2} \psi_j$ and equating like powers of $\sqrt\eps$ in $\rPm\psi=0$ provides the following constraints:
\begin{align}\label{eq:asympt}
  \rPmz \psi _0 =0,\quad \rPmz\psi_1+\rPmo\psi_0=0, \quad \rPmz\psi_{j+2} +\rPmo \psi_{j+1} + \rPmt\psi_j =0, \ j\geq0.
\end{align}
We do not explicitly keep track of the dependence of $\psi_j$ on $m$. To characterize the operators $\rPmj$ and solutions to the above equations, we first introduce an appropriate functional setting.
\paragraph{Functional setting.} 
%

We need to construct a number of spaces of functions of $z$, functions of $(t,x,\xi)$, and functions of $(t,x,\xi,z)$.  
Let $p\in\Nm$ and $\fa=\partial_z+z$ with $\fa^*=-\partial_z+z$. For functions $\psi(z)$, we define the spaces $\mS_p$ with the norm
\begin{align*}
  \|\psi\|_{p}= \|(\fa^*)^{p} \psi\|_0,\qquad \|\cdot\|_0= \|\cdot \|_{L^2(\Rm)}.
\end{align*}
These are the standard Hilbert spaces associated with decomposition into a basis of Hermite functions. For functions $\psi(t,x,\xi)$, we define the spaces $\mS_p$ with the norm 
\begin{align*}
  \|\psi\|_{p}= \sup_{t\in[-\rT,\rT]} \sum_{\sum_{j=1}^3\alpha_j=p} \| D_x^{\alpha_1}D_t^{\alpha_2} \aver{\xi}^{\alpha_3} \psi\|_0,\qquad \|\cdot\|_0= \|\cdot \|_{L^2(\Rm^2)},\qquad \aver{\xi}:=\sqrt{1+\xi^2}.
\end{align*}

Since the homogeneities of the operators that appear in the theory are different in the $z$ and $\xi$ variables (all functions depending smoothly on $(t,x)$ at the multi-scale level), we define for $p=(p_1,p_2)$ with $p_j\in\Nm$ for $j=1,2$, the spaces $\mCS_p$ with the norm
\begin{align*}
  \|\psi\|_{p}= \sup_{t\in[-\rT,\rT]} \sum_{\sum_{j=1}^4\alpha_j=p_2} \|(\fa^*)^{p_1+\alpha_1} D_x^{\alpha_2}D_t^{\alpha_3} \aver{\xi}^{\alpha_4} \psi\|_0,\qquad \|\cdot\|_0=\|\cdot\|_{L^2(\Rm^3)}.
\end{align*}

All spaces $\mS_p(\Rm^q;\Cm^l)$ and $\mCS_p(\Rm^q;\Cm^l)$ of vector-valued functions $\psi$ are defined similarly component-wise. 
The norm of an operator $A$ from $\mS_{p'}$ to $\mS_{p}$ or from $\mCS_{p'}$ to $\mCS_p$ is denoted by 
\begin{align*}
   \vvvert A \vvvert_{p,p'}.
\end{align*}
We denote by $\mS_p$ or $\mCS_p$ the corresponding spaces in all above cases when the context is clear.

\paragraph{Operator assumptions.}
We now collect the assumptions we need to make on ($L$ via) $\rPmj$.
\begin{enumerate}
\item[(i)] The operator $P_0$ in \eqref{eq:P0} is a polynomial of degree $q$ in the last four variables.
\item[(ii)] We assume the existence of a real-valued solution $G_m(t,x,\xi;x_0)$ of \eqref{eq:eikonal} of the form
\begin{align}\label{eq:Gmglobal}
   G_m(t,x,\xi;x_0) = -B_m(\xi;x_0) t +A_m(x,\xi;x_0)
\end{align}
with $A_m$ and $B_m$ smooth and uniformly Lipschitz in the $\xi$ variable for $\xi\in\Xi=\Xi(m)$ an open subset of $\Rm$, and such that $\partial_x A_m(x_0,\xi;x_0)=\xi$. 

Associated to these branches are the one-dimensional eigenspace 
\begin{align}\label{eq:kernelN}
 N= N_{x,\xi,m;x_0}:=\Cm \phi_m(x,\xi,z;x_0)\subset L^2(\Rm_z;\Cm^l)
\end{align}
with $\phi_m(x,\xi,z;x_0)$ in $\mCS_p$ for any $p=(p_1,p_2)\in\Nm^2$, uniformly in $x_0$.  Solutions of $\rPmz \psi_0=0$ in $\mCS_p$ for $p=(p_1,p_2)$ are given by $\psi_0(t,x,\xi,z;x_0) = f(t,x,\xi;x_0) \phi_m(x,\xi,z;x_0)$ for $f\phi_m\in\mCS_{p}$ arbitrary, or equivalently $f\in \mS_{p_2}$ arbitrary.

\item[(iii)] For an integer $q_2\geq0$, the equation $\rPmz \psi=g$ admits solutions  in $\mCS_p$ when $g\in N^\perp \cap \mCS_{p+q_2}$ using the orthogonal decomposition $L^2(\Rm;\Cm^l)=N\oplus N^\perp$. When $g\in N^\perp \cap \mCS_{p+q_2}$, the unique solution in $N^\perp  \cap \mCS_{p}$ is denoted by $\rPmzi g$. General solutions in $\mCS_{p}$ are given by $\rPmzi g + f_1(x,\xi;x_0)\phi_m(x,\xi,z;x_0)$ for an arbitrary $f_1$ such that $f_1\phi_m\in\mCS_p$. 

The operator $\rPmzi$ is bounded  from $N^\perp \cap \mCS_{(p_1,p_2+q_2)}$ to $N^\perp \cap \mCS_{(p_1,p_2)}$ by $C_p$ independent of $m$:
\begin{align}\label{eq:boundrPmzi}
   \vvvert \rPmzi \vvvert_{(p_1,p_2),(p_1,p_2+q_2)} \leq C_p.
\end{align}

\item[(iv)] For $g(t,x,\xi,z;x_0)=-\rPmo( f(t,x,\xi;x_0)\phi_m(x,\xi,z;x_0))$, then $g\in N^\perp$ is equivalent to
\begin{align}\label{eq:mTm}
   \mT_m f := (\phi_m,\rPmo (f\phi_m))_z =0
\end{align}
where $(\cdot,\cdot)_z$ is the standard inner product on $L^2(\Rm_z;\Cm^l)$. The equation 
\begin{align*}
  \mT_m f=g,\qquad f(0,x,\xi;x_0)=h(x,\xi;x_0)
\end{align*}
with $h\in\mS_p$ and $g\in \mS_p$ admits a unique solution $f(t,x,\xi;x_0)$ in $\mS_p$ with a constant $C_p$ independent of $m$ and $\rT$ such that for each $p=(p_1,p_2)$ we have
\begin{align}\label{eq:stabtr}
   \|f\phi_m\|_p  \leq C_p (\|h\phi_m\|_p + \aver{\rT} \|g\phi_m\|_p),\qquad \aver{\rT}=\sqrt{1+\rT^2}.
\end{align}
Moreover, we assume the following finite speed of propagation result. There is a universal constant $C$ such that if $g$ and $h$ are supported in an interval $[x_1,x_2]$ in the variable $x$, then $f(t,x,\xi;x_0)$ is supported in the interval $[x_1-C|t|,x_2+C|t|]$ in the $x$ variable.

\item[(v)] The operators $\rPmj$ are bounded uniformly in $m$ from $\mCS_{(p_1+j,p_2+q)}$ to $\mCS_p$, i.e.,
\begin{align}\label{eq:boundrPmj}
   \vvvert \rPmj \vvvert_{(p_1,p_2),(p_1+j,p_2+q)} \leq C_p.
\end{align}
\end{enumerate}
\begin{remark} \label{rem:integers} \rm
 For the Dirac problem, we have $q=q_2=1$ and $l=2$ while for the Klein-Gordon problem, we have $q=2$, $q_2=0$, and $l=1$.
\end{remark}
\begin{remark} \label{rem:generalties} \rm
  The asymptotic construction of wavepackets is based on a classical two-scale separation with microscopic oscillation parametrized by $z$ away from the interface and $(\xi,m)$ along the interface, and macroscopic oscillations parametrized by time $t$ and $x$ along the interface. The leading eikonal equation \eqref{eq:eikonal} reflects this combination of confinement away from the interface and propagation along it in the microscopic variables. Hypotheses (i)-(ii) allow us to construct solutions to this eikonal equation. The behavior in the macroscopic variables $(t,x)$ is obtained by the next-order term $\rPmo$. Because of confinement, macroscopic oscillations are occurring purely along the interface and are described by the operator $\mT_m$ in \eqref{eq:mTm} integrating out the transverse variable $z$. Hypotheses (iii)-(iv) ensure that $\rPmz$ is appropriately invertible so that higher-order terms in the expansion may be computed. Finally, hypothesis (v) ensures that errors resulting from the above approximation are controlled in a square integrable sense. 
\end{remark}
\subsection{Wavepacket construction}
The construction is based on solving the equations \eqref{eq:asympt} in turn. Let $J$ be a fixed integer quantifying the accuracy of the wavepacket. 

The leading order equation $\rPmz \psi _0 =0$ is solved based on (ii) above and provides for each branch $m$ the family of solutions 
\begin{align*}
  \psi_0(t,x,\xi,z;x_0) = f_0(t,x,\xi;x_0) \phi_m(t,x,\xi,z;x_0),
\end{align*}
where $f_0\in \mCS_p$ is arbitrary at this stage. We do not keep the dependence in $m$ explicit.  

The next-order equation $\rPmz\psi_1+\rPmo\psi_0=0$ provides both an equation to solve for $\psi_1$ up to a term $f_1(t,x,\xi;x_0) \phi_m(t,x,\xi,z;x_0)$ as indicated in (iii) as well as an equation to solve for $f_0$ in terms of the compatibility condition $\rPmo(f_0\phi_m)\in N^\perp$.

Let 
\begin{equation}\label{eq:initcond}
  \hat f_m(x,\xi;x_0) \quad \mbox{be  such that } \quad \hat f_m\phi_m\in \mCS_{p} \quad \mbox{ is a {\em given} initial condition}
\end{equation}
with support in the $\xi$ variable in $\Xi$ introduced in (ii). It is the only free term in the forthcoming construction of the wavepackets. 

We define $f_0$ as the solution to $\mT_mf_0=0$ with initial condition $f_0(t=0)=\hat f_m$. We deduce from (iv) that $f_0\phi_m\in\mCS_{p}$ with $\|f_0\phi_m\|_p\leq C \|\hat f_m\phi_m\|_p$. 

Once $f_0$ is constructed, we solve for $\rPmz\psi_1+\rPmo\psi_0=0$ to obtain using (iv) and (v) that $\psi_1=\tilde\psi_1  + f_1\phi_m$ with $\tilde\psi_1 =-\rPmzi \rPmo \psi_0 \in \mCS_{(p_1-1,p_2-q-1)}$ and $f_1$ arbitrary. The next equation is  $\rPmz\psi_2+\rPmo\psi_1+\rPmt\psi_0=0$, from which we deduce the compatibility condition
\begin{align*}
  \mT_m f_1 = -(\phi_m,\rPmo\tilde\psi_1 + \rPmt \psi_0)_z.
\end{align*}
Augmented with $f_1(t=0)=0$ and using \eqref{eq:stabtr} and \eqref{eq:boundrPmj}, it admits a unique solution such that:
\begin{align*}
   \|f_1\phi_m\|_p \leq C_p \aver{\rT} \|\rPmo\tilde\psi_1 + \rPmt \psi_0\|_p \leq C_p  \aver{\rT}  \|\hat f_m\phi_m\|_{(p_1-2,p_2-(2q+q_2))}.
\end{align*}

We construct the higher-order terms $\psi_j=\tilde\psi_j+f_j\phi_m$ iteratively from 
 $\rPmz \psi_{j+1} + \rPmo\psi_j + \rPmt \psi_{j-1}=0$.
Assume $\psi_j$ constructed with $f_j$ defined such that the compatibility condition $(\phi_m,\rPmo\psi_j + \rPmt \psi_{j-1})_z=0$ holds and that moreover,
\begin{align}\label{eq:regpsij}
   \|\tilde\psi_j\|_p & \leq C \aver{\rT}^{j-1} \|\hat f_m\phi_m\|_{p_1-2j+1,p_2-(2q+q_2)j+q} \\ \nonumber
  \|\psi_j\|_p +  \|f_j\phi_m\|_p &\leq C  \aver{\rT}^{j} \|\hat f_m\phi_m\|_{p_1-2j+1,p_2-(2q+q_2)j+q} .
\end{align}
These bounds have been verified above for $j=1$. 
Then we construct in turn
\begin{align*}
   \tilde\psi_{j+1} &= -\rPmzi (\rPmo\psi_{j} + \rPmt\psi_{j-1})\quad \mbox{ and } \quad 
   \mT_m f_{j+1} = -(\phi_m, \rPmo\tilde\psi_{j+1} + \rPmt \psi_{j})_z
\end{align*}
with the last equation augmented with initial conditions $f_{j+1}(t=0)=0$.

Using \eqref{eq:regpsij} and \eqref{eq:boundrPmzi} as well as \eqref{eq:boundrPmj}, we obtain the bound for $\tilde\psi_j$ given in \eqref{eq:regpsij} with $j$ replaced by $j+1$. Using \eqref{eq:boundrPmj} as well as  \eqref{eq:stabtr}, we next obtain the bound for $f_{j+1}\phi_m$ given in \eqref{eq:regpsij}.

This allows us to construct $\psi_j$ for $0\leq j\leq J$ to arbitrary order $J$ provided that $p_1$ and $p_2$ are sufficiently large. We conclude the construction by setting $f_J=0$ for concreteness so that $\psi_J=\tilde\psi_J$ has the same regularity as $\tilde\psi_J$. In particular,
\begin{align*}
   \|\psi_J \|_{p_1,p_2} \leq C \aver{\rT}^{J-1} \|\hat f_m\phi_m \|_{p_1+2J-1,p_2+(2q+q_2)J-q}.
\end{align*}

We define the complete wavepacket as 
\begin{align}\label{eq:psimJ}
  \psi_m^J (t,x,\xi,z;x_0)= \sum_{j=0}^J \eps^{\frac j2} \psi_j (t,x,\xi,z;x_0).
\end{align}

Note that $\psi_m^J$ has the same regularity as $\psi_J$ although it is better to see it as a superposition of terms whose regularity decreases with their amplitude. Note also that 
\begin{align*}
  \psi_m^J (0,x,\xi,z;x_0)=  \hat f(x,\xi;x_0) \phi_m(x,\xi,z;x_0) +\sum_{j=1}^J \eps^{\frac j2} \tilde \psi_j (0,x,\xi,z;x_0).
\end{align*}
The initial condition of the wavepacket $\psi_m^J$ is therefore given only to leading order by the prescribed $\hat f(x,\xi;x_0) \phi_m(x,\xi,z;x_0)$ in \eqref{eq:initcond}. The remainder involving the above sum over $\tilde\psi_j$ is necessary to ensure that the initial wavepacket $\psi_m^J(0)$ belongs to the perturbed $m$th branch with sufficient accuracy and {\em not} the unperturbed branch parametrized by $\phi_m$. The terms $\tilde \psi_j (0,x,\xi,z;x_0)$ decompose in $L^2(\Rm_z;\Cm^l)$ over the eigenfunctions $\phi_{m'}$ with $m'\not=m$.

We now construct $\Psi^J$ from $\psi^J$ using \eqref{eq:Psi} and observe that 
\begin{align} \nonumber
  \tilde P_\eta \Psi^J(t,x,z) & =\eps^{-\frac14}\dint_{\Xi} e^{\frac i{\sqrt\eps} G_m} \rPm\psi^J (t,x,\xi,z) d\xi  = \eps^{-\frac14} \dint_{\Xi} e^{\frac i{\sqrt\eps} G_m}\eps^{\frac{J+1+q}2} \rr_J (t,x,\xi,z) d\xi  \\  & =: \eps^{\frac{J+1+q}2} r_J(t,x,z;x_0) , \label{eq:rj}
\end{align}
where the multiscale corrector $\rr_J$ satisfies
\begin{align}\label{eq:rrJ}
 \rr_J = -(\rPmo+\sqrt\eps \rPmt)\psi_J - \rPmt \psi_{J-1}.
\end{align}
Therefore, using \eqref{eq:boundrPmj} and \eqref{eq:psimJ}, we derive the bound
\begin{align}\label{eq:bdrrJ}
  \|\aver{\xi} \rr_J\| \leq C  \aver{\rT}^{J-1}\big( \|\hat f_m\phi_m\|_{2J,(2q+q_2)J+1} + \sqrt\eps  \|\hat f_m\phi_m\|_{2J+1,(2q+q_2)J+1} \big)
\end{align}
with here $\|\cdot\|$ the norm with $p=(0,0)$. We thus obtain the following result:
\begin{lemma}\label{lem:psimJ}
   Let $\hat f_m(x,\xi;x_0)$ be given such that  $\hat f_m\phi_m\in \mCS_p$ in \eqref{eq:initcond} for $p_1\geq 2J+1$ and $p_2\geq (2q+q_2)J$. We also assume that $\hat f_m$ vanishes when $|x-x_0|\geq\eta$. 
   The wavepacket $\psi_m^J$ is constructed as in \eqref{eq:psimJ} with terms $\psi_j$ satisfying \eqref{eq:regpsij}. 
   Let $r_J$ be defined in \eqref{eq:rj}. 
   Then, we have the estimate:
\begin{align*}
  \|r_J\|:=  \sup_{t\in [-\rT,\rT]} \|r_J(t,\cdot)\|_{L^2(\Rm^2;\Cm^2)} \leq C   \eps^{-\frac14}\aver{\rT}^{J-1} \|\hat f_m\phi_m\|_{2J+1,(2q+q_2)J}.
\end{align*}
\end{lemma}
\begin{proof}
This uses that $G_m$ is real-valued and that $\rr_J=\aver{\xi}^{-1} \aver{\xi}r_J$ is integrable by Cauchy-Schwarz and \eqref{eq:bdrrJ}.
\end{proof}
The above estimate in $\eps^{-\frac14}$ is not optimal and may be replaced by an estimate independent of $\eps$ with a more detailed analysis of the oscillatory properties of $G_m$. Since $J$ may be chosen arbitrarily large modulo appropriate regularity assumptions on $\hat f_m$, we use  this convenient estimate.

The initial condition is assumed to have (small) compact support and (square) integrability in $x$. Since the initial conditions we are interested in are supported in the $\sqrt\eps-$vicinity of $(x_0,0)\in\Phi^{-1}(\Gamma)$, there is no loss in multiplying it by $\chi_\eta(x-x_0)$, say.

\subsection{Wavepacket Analysis}
\label{sec:wpa}
%


We next define
\begin{align*}
  v^J(t,x,y)=\eps^{-\frac14}\Psi^J(t,x,z)
\end{align*} 
in the rectified variables $(x,y)$ dropping the dependence in $m$ to simplify. It remains to construct a wave-packet in the original variables $\Phi(x,y)$. When $\Gamma$ is open and $\Phi$ is a global diffeomorphism, this is achieved by constructing $u^J=\Phi^{-*}v^J=v^J\circ \Phi^{-1}$. When $\Gamma$ is closed, however, $\Phi^{-1}$ is defined globally only if we identify any point $(x,y)$ with $(x+k\rL,y)$ in the rectified variables. We thus introduce the periodic wavepacket
\begin{align*}
  v^J_\sharp(t,x,y) = \dsum_{k\in\Zm} v^J(t,x+k\rL,y).
\end{align*}
More generally, we denote by $v_\sharp$ the periodization of $v$ as given above.
Then $\Phi$ is a diffeomorphism from $[0,\rL)\times [-2\eta,2\eta]$ to its image and we then define
\begin{align}\label{eq:ufromv}
  u^J=\Phi^{-*}v^J_\sharp
\end{align} 
in the original variables. To simplify notation, we also denote by $v^J_\sharp$ the function $v^J$ when $\Gamma$ is open so that the above expression holds both when $\Gamma$ is open and when it is closed. Thanks to the presence of $\chi_\eta$ in the definition of $P_\eta$, the function $v^J$ is supported in the $2\eta-$vicinity of the axis $y=0$ while $u^J$ is also supported in an $\eta-$vicinity of $\Gamma$ after applying the smooth (inverse) pullback.
\subsubsection{Accuracy of wave-packet parametrix.}

Let now $u$ be the solution of $Lu=0$ with initial conditions $(\sqrt\eps\partial_t)^j u_{|t=0}=(\sqrt\eps\partial_t)^ju^J_{|t=0}$ for $0\leq j\leq q-1$. Then
\begin{align}\label{eq:tilderJ}
 \eps^{-\frac q2} L (u-u^J) = -  \eps^{-\frac q2}L u^J = -\eps^{\frac {J+1-q}2} \Phi^{-*}(\tilde r_J)\ \ \mbox{ for } \ \ \tilde r_J(t,x,y):=\sum_{k\in\Zm} \eps^{-\frac14}r_J(t,x+k\rL,\frac y{\sqrt\eps})
\end{align}
and the convention that the latter sum equals $ \eps^{-\frac14}r_J(t,x,\frac y{\sqrt\eps})$ when $\rL=\infty$, i.e., when $\Gamma$ is open.

We make the following assumption on the operator $L$.
\begin{enumerate}
    \item[(vi)] The equation $Lu=g$ augmented with initial conditions $(\sqrt\eps \partial_t)^ju_{|t=0}=0$ for $0\leq j\leq q-1$ admits a unique solution with the following energy stability estimate:
    \begin{align}\label{eq:controlmL}
         \mL[u](t) \leq C \eps^{\frac12-q} \dint_0^t \|g(s,\cdot)\|_{L^2(\Rm^2)} ds.
    \end{align}
\end{enumerate}

For the Dirac equation with $q=1$, the above (square root of the) energy functional is nothing but the $L^2$ norm. For Klein-Gordon equations with $q=2$, the above functional is a standard (semiclassical) energy functional associated with second-order wave equations. In both cases, the energy functional controls the following quantities
    \begin{align*}
       \|(\sqrt\eps \partial_t)^{q-1} u(t,\cdot) \|_{L^2(\Rm^2;\Cm^l)} \leq \mL[u](t) .
    \end{align*}

\begin{theorem} \label{thm:localerror}
  Under assumptions (i)-(vi) and for wavepackets based on an initial condition $\hat f_m\in \mCS_p$ for $p_1\geq 2J+1$ and $p_2\geq (2q+q_2)J$ with support in $\xi$ in $\Xi$ introduced in (ii), let $u$ be the solution of $Lu=0$ with initial conditions $(\sqrt\eps\partial_t)^j u(t=0)=(\sqrt\eps\partial_t)^ju^J(t=0)$ for $0\leq j\leq q-1$. Then
  \begin{align}\label{eq:estimerror}
     \mL[u-u^J](t) \leq C   \mE(\eps,J,\rT)  \|\hat f_m\phi_m\|_p,\qquad  \mE(\eps,J,\rT) :=\eps^{-\frac14} \aver{\rT}^c  \eps^{1-q}  \big(\aver{\rT}\eps^{\frac12}\big)^J,
  \end{align}
  provided $p$ is as in Lemma \ref{lem:psimJ}.
  Here, $c=0$ when $\Gamma$ is open while $c=1$ when $\Gamma$ is closed.
\end{theorem}
\begin{proof}
 The number of terms in $\tilde r_J$ defined in \eqref{eq:tilderJ} grows at most linearly in time by finite speed of propagation of the transport equation. This shows that $\tilde r_J$ in the $L^2$ sense is at most $\aver{\rT}^c$ larger than $r_J$ in the same sense. It then remains to use (vii) and the bound in Lemma \ref{lem:psimJ}, as well as the fact that $\Phi^{-*}$ is a smooth diffeomorphism to conclude.
\end{proof}
Since $J$ above is arbitrary, at least when the initial condition is sufficiently smooth, we obtain a control on the wavepackets (morally) up to times such that $\eps^{\frac12}\aver{\rT}\ll1$.

\subsubsection{Arbitrary localized initial conditions.}

The above theorem shows that for appropriate choices of initial conditions \eqref{eq:initcond} in a given $m$th branch, then the wavepacket $u^J=u^J_m$ is an accurate approximation of the solution to the wave equation $Lu=0$ with the same initial conditions. The ansatz \eqref{eq:Psi} also offers a valuable tool to obtain qualitative properties of the wavepacket as we will demonstrate below.

Somewhat hidden in the assumptions of the above result is that the phase functions $G_m$ are well-behaved only for $\xi\in\Xi$. This restricts the choice of $\hat f_m(x,\xi;x_0)$. Another constraint is the explicit form of the initial condition for $u^J$.  In the rectified coordinates and using (ii), we find
\begin{align*}
  v^J(0,x,y) = \frac{1}{\sqrt\eps} \dint_{\Xi} e^{\frac i{\sqrt\eps} A_m(x,\xi;x_0)} \psi^J_m(0,x,\xi,\frac y{\sqrt\eps}) d\xi.
\end{align*}
The above phase satisfies $\partial_x A_m(x_0,\xi;x_0)=\xi$ so that $A_m$ is asymptotically close to $\xi(x-x_0)$ but not exactly that function. Thus, $\psi_m^J$ resembles the semiclassical Fourier transform of $v^J$ but is not exactly defined as such.

In favorable settings still of practical interest, we may choose $\Xi=\Rm$ (or $\Xi=\{|\xi|>\delta\}$ for arbitrary $\delta>0$)  in which case (essentially) no restriction beyond regularity is imposed on $\hat f_m$, and moreover, $A_m$ may be transformed to $k(x-x_0)$ by a smooth global change of variables. In such a setting, the solution of $Lu=0$ with (semiclassical) general initial conditions $u(0,x,y)$ localized in the vicinity of $(x_0,y_0)\in\Gamma$ may be represented as a superposition of wavepackets of the form $\psi^J_m$. Obtaining such a result is the main objective of this section.

Besides (i)-(vi), the main hypotheses we make are the following:
\begin{enumerate}
  \item[(vii)]   We assume the existence of a countable number $m\in M$ such that (i)-(v) holds. At each $(x,\xi;x_0)$ fixed, $M=\cup_{j=1}^q M_j$ may be decomposed as a disjoint union of $q$ sets $M_j$ for $1\leq j\leq q$ such that for each $j$, $\{\phi_m\}_{m\in M_j}$ is an orthonormal basis of $L^2(\Rm_z;\Cm^l)$ independent of $j$.  
  
  The open set introduced in (ii) is given by $\Xi=\{|\xi|> \xi_0\}$ for some $\xi_0>0$ or $\Xi=\Rm$ independent of $m$. Moreover, for each $m$, we have a change of variables $(x,\xi)\to (x,k)$ with
\begin{align*}
       k=\frac{A_m(x,\xi;x_0)}{x-x_0}
\end{align*} 
that is a regular diffeomorphism from $\Rm\times\Xi$ to its image that includes a domain of the form $\Rm\times K$ with $K=\{|k|>k_0\}$ for $k_0>0$ or $K=\Rm$ for a domain $K$ independent of $m$. 
 \item[(viii)] Let $M \ni m=(n,\ell)\equiv n\ell$ be the decomposition of the branch indices in (vii) and $B_m\equiv B_{n\ell}$. For each fixed $n$, let $\rA_{\ell j}(\xi;x_0)=(-B_{n\ell}(\xi;x_0))^{(j-1)}$ be a $q\times q$ Vandermonde matrix. We assume that $\rA$ is uniformly invertible for $\xi\in\Xi$.
\end{enumerate}
In the cases we will consider in practice, $A_m$ and the change of variables are in fact independent of $m$. We will be able to choose $\Xi=K=\Rm$ for the Dirac operator and $\Xi=\{|\xi|>\delta\}$ for arbitrary $\delta>0$ for the Klein-Gordon operator, in which case $k_0$ above may be chosen arbitrarily small as well. We then have the following result:
\begin{theorem}\label{thm:globalerror}
  Assume that hypotheses (i)-(viii) hold. Let $(x_0,0)\in \Phi^{-1}(\Gamma)$ and consider initial conditions $(\sqrt\eps\partial_t)^ju_{|t=0}(x,y)$ concentrated in the vicinity of $\Phi(x_0,0)$ and represented in the rectified coordinates $(x,z)$ by
   \begin{align*}
      (\sqrt\eps D_t)^j \Psi(0,x,z) = \eps^{-\frac14} \dint_{K} e^{\frac i{\sqrt\eps} k(x-x_0)} \hat \Psi_j(k,z) dk,\qquad 0\leq j\leq q-1,
   \end{align*}
   with $\hat \Psi_j(k,z)$ thus the semiclassical Fourier transforms of the initial conditions in the first variable. We assume that  $(\sqrt\eps D_t)^j \Psi(0,x,z)\in \mCS_{p}$ for $(p_1,p_2)$ sufficiently large and that they are supported in the $\eta-$vicinity of $x_0$ in the variable $x$.
   
With $u=\Phi^{-*}v$ and $v(t,x,y)=\eps^{-\frac14}\Psi(t,x,\frac{y}{\sqrt\eps})$, the above provides arbitrary (sufficiently smooth, localized and with frequency content in $K$) initial conditions for $u$ solutions of $Lu=0$.
   
Then we may construct a superposition of wavepackets 
\begin{align*}
      u^J(t,x,y) = \Phi^{-*} \dsum_m v_m^J,\qquad v_m^J(t,x,y)= \eps^{\frac14} \Psi_m^J(t,x,\frac{y}{\sqrt\eps}),
\end{align*}
where $\Psi_m^J$ is constructed from $\psi_m^J$ using \eqref{eq:Psi} and $\psi_m^J$ is constructed based on appropriately chosen ($\eps-$dependent) $\hat f_m\in \mCS_p$ as described in Lemma \ref{lem:psimJ}. Moreover, with $\mE$ as in \eqref{eq:estimerror}, we have:
   \begin{align}\label{eq:arbiniterror}
     \mL[u-u^J](t) \leq C \mE(\eps,J,\rT).
  \end{align}
\end{theorem}
The construction of the terms $\hat f_m$ from knowledge of $(\eps D_t)^j \Psi(0,x,z)$ is explicit but complicated. This theorem shows that {\em any sufficiently smooth initial condition} localized in the vicinity of a point $(x_0,y_0)\in\Gamma$ (in the un-rectified coordinates) may be decomposed over wavepackets whose dynamics are predicted by an ansatz of the form \eqref{eq:Psi}.

The constant $p=(p_1,p_2)$ in the above theorem is larger but not fundamentally different from that given in Theorem \ref{thm:localerror}. Since its expression does not seem particularly meaningful, we do not attempt to obtain it explicitly as this shortens the derivation.


\begin{proof}
The proof is based on two main ingredients. The first one is a change of variables from $A(x,\xi;x_0)$ to $k(x-x_0)$ to show that the initial condition for the wavepackets of Lemma \ref{lem:psimJ} may be transformed to the semiclassical Fourier transform of arbitrary (sufficiently smooth) initial conditions. Second, we observe that the constructed wavepackets satisfy prescribed initial conditions only up to an error of order $O(\sqrt\eps)$. We then iterate the above procedure to construct additional wavepackets that provide an overall error of order $O(\eps^{\frac J2})$.

Introduce the ansatz for sufficiently smooth functions $f_m$:
\begin{align}\label{eq:tildePsi}
  \tilde \Psi(t,x,z) &=  \eps^{-\frac14} \dint_{\Xi} \sum_{m\in M} e^{\frac i{\sqrt\eps} G_m(t,x,\xi)} f_m(t,x,\xi)\phi_m(x,\xi,z) d\xi \\ \nonumber
    (\sqrt\eps D_t)^j \tilde\Psi (0,x,z) &=   \eps^{-\frac14} \dint_{\Xi} \sum_{m\in M} e^{\frac i{\sqrt\eps} A_m(x,\xi)} (\partial_t G_m)^j  f_m(0,x,\xi)\phi_m(x,\xi,z) d\xi + O(\sqrt\eps)
\end{align}
where the $O(\sqrt\eps)$ term in $\mCS_p$ spaces involves derivatives in time of $f_m(t,x,\xi)$ evaluated at $t=0$. 

Our first objective is to find $f_m(0,x,\xi)$ such that for $0\leq j\leq q-1$,
\begin{align*}
 \dint_{\Xi} \sum_{m\in M} e^{\frac i{\sqrt\eps} A_m(x,\xi)} (-B_m(\xi))^j  f_m(0,x,\xi)\phi_m(x,\xi,z) d\xi  =  \dint_{K} e^{\frac i{\sqrt\eps} k(x-x_0)} \hat \Psi_j(k,z) dk.
\end{align*} 
We write $m=(n,\ell)\equiv n\ell$ with $\phi_m=\phi_n$ by assumption (vii) an orthonormal basis of $L^2(\Rm_z;\Cm^l)$ and $1\leq \ell\leq q$, and (viii) to recast the above left-hand side as 
\begin{align*}
 \dsum_n \dint_{K} e^{\frac i{\sqrt{\eps}}k(x-x_0)} \Big(  \dsum_{\ell}(-B_m(\xi))^j f_m(0,x,\xi) \partial_k\xi \Big) \phi_n(x,\xi,z) dk.
\end{align*}
Above, we implicitly use that $\xi=\xi(x,k)$ and Jacobian $0<\partial_k\xi=\partial_k\xi(x,k)$.
Moreover, since $\phi_n(x,\xi,\cdot)$ is an orthonormal basis of $L^2(\Rm_z;\Cm^l)$,
\begin{align*}
   \hat \Psi_j(k,z) = \dsum_{n} \tilde f_{nj}(k) \phi_n(x,\xi(x,k),z).
\end{align*}
Therefore, we solve our problem provided 
\begin{align}\label{eq:VDM}
  \dsum_\ell  (-B_{n\ell}(\xi))^j f_{n\ell} (0,x,\xi)  \partial_k\xi=  \tilde f_{nj}(k).
\end{align}
This is Vandermonde $q\times q$ system (as in the construction of parametrices for wave equations; see, e.g., \cite{GS-CUP-94}) that admits a unique solution for each $k\in K$ by assumption (viii) and definition of $K$. When $q=1$, we note that this is a scalar equation with explicit solution:
\begin{align*}
   f_m(0,x,\xi) = \partial_\xi k (x,\xi) \tilde f_m(x,\xi).
\end{align*}
Note that since $\Psi(0,x,z)$ is supported in the $\eta-$vicinity of $x_0$, we also have that $f_m(0,x,\xi)$ is supported in the vicinity of $x_0$.

We have thus constructed an ansatz $\tilde \Psi(t,x,\xi)$, which at time $t=0$ is a superposition over $m$ of initial conditions $\hat f_m(x,\xi)\phi_m(x,\xi,z)$ that allow one to construct wavepackets $\psi_m^J$ as summarized in Lemma \ref{lem:psimJ}. The functions $f_m(t,x,\xi)$ are defined as in the construction of the wavepackets. 

Each wavepacket $\psi_m^J$ provides an error in the norm $\mL$ of order $\mE \|\hat f_m\phi_m\|_p$ as provided in Theorem \ref{thm:localerror}. The total error in the same sense is therefore of order $\mE$ provided that 
\begin{align*}
  \dsum_m  \|\hat f_m\phi_m\|_p < \infty.
\end{align*}
This condition is satisfied when $\hat f_m(\xi)$ decays sufficiently rapidly in $m$ and $\xi$, which in turn is clear when $\hat\Psi_j(k,z)$ is sufficiently rapidly decaying in $k$ and sufficiently smooth and rapidly decaying in $z$. This holds since $p=(p_1,p_2)$ is assumed sufficiently large.

So, for $u^J_1$ constructed with the above $\Psi_m^J$, we obtain \eqref{eq:arbiniterror}. The initial conditions for $u^J_1$ differ from those prescribed for $u$ in two ways.  First, the $\sqrt\eps-$ error term in \eqref{eq:tildePsi} involves time derivatives of $f_m$ at time $t=0$. Since $f_m$ may be chosen arbitrarily regular, the sum over $m$ of these contributions are bounded in the energy sense $\mL[u]$ at time $0$ once pulled back to the unrectified coordinates; to obtain this, we use the bound on $\mL$ prescribed in (vi). We thus obtain an error in the energy sense between the initial conditions for $\Psi$ and those for $\tilde\Psi$ of order $\sqrt\eps$. Second, the wavepacket $u^J_1$ satisfies the initial conditions prescribed by $\tilde\Psi$ also up to an error of order $\sqrt\eps$. 

\medskip

We have thus constructed $\Psi^J_1$ with initial conditions  $(\sqrt\eps D_t)^j\Psi^J_1(0,x,z)= (\sqrt\eps D_t)^j\Psi(0,x,z)+\sqrt\eps \Psi_{2j}(x,z)$ with $\Psi_{1j}$ in any $\mS_{p'}$ provided that $p$ is sufficiently large. We may therefore define $\hat \Psi_{2j}(k,z)$  as the semiclassical Fourier transform of  $\Psi_{2j}(x,z)$. It remains to repeat the above procedure to construct initial conditions $\hat f_{2,m}(\xi,x)$ and corresponding wavepackets as prescribed in Lemma \ref{lem:psimJ} and $u^{J-1}_2$ mimicking the construction of $u^J_1$. Iterating this process $J-1$ times provides a superposition of wavepackets 
\begin{align*}
 u =  u^J_1 + \sqrt\eps u^{J-1}_2 + \ldots \eps^{\frac{J-1}2} u^1_{J} + O(\eps^{\frac J2})
\end{align*}
with an error in the initial condition of order $\eps^{\frac J2}$ in the $\mL$ sense.  Applying \eqref{eq:arbiniterror} as we have done above for the term $u^J_1$, each of these wavepackets provides an error at time $t$ of order 
\begin{align*}
   \eps^{\frac{J-1}2} \mE(\eps,J-j,\rT) \leq \mE(\eps,J,\rT).
\end{align*}
This concludes the derivation of the theorem.
\end{proof}

In the favorable setting of (vii)-(viii) where $G_m$ may be chosen linear in $t$ and $G_m(t=0)$ may be bijectively transformed to the phase of a standard semiclassical Fourier transform, we thus observe that a superposition over $m$ of terms of the form \eqref{eq:Psi} allows for the propagation of arbitrary (sufficiently smooth) initial conditions localized in the vicinity of a point $(x_0,y_0)\in \Gamma$, up to an arbitrarily small error in powers of the semiclassical parameter $\eps$.  See remark \ref{rem:turning} below for brief comments on the setting where (vii)-(viii) does not hold.


%
%
%
\subsubsection{Leading term in relativistic wavepackets.}

Consider the setting of Theorem \ref{thm:localerror}. The wavepacket $u^J=\Phi^{-*} v_\sharp^J$ with $v^J(t,x,y)=\eps^{-\frac14} \Psi^J(t,x,z)$ and $y=\sqrt\eps z$ admits a leading term $u_0=\Phi^{-*} v_{\sharp 0}$ with $v_0(t,x,y)= \eps^{-\frac14} \Psi_0(t,x,z)$ for 
\begin{align} \label{eq:Psi0}
 \Psi_0(t,x,z ) =  \eps^{-\frac14} \dint_{\Xi} e^{\frac i{\sqrt\eps} G_m(t,x,\xi)} f_0(t,x,\xi;x_0) \phi_m (x,\xi,z;x_0) d\xi.
\end{align}
Moreover, from Theorem \ref{thm:localerror} and (vii), we have the estimate
\begin{align}\label{eq:erroru0}
  \mL[u-u_0]+\mL[u^J-u_0]\leq C\sqrt\eps.
\end{align}

Non-dispersive modes are characterized by a phase function of the form
\begin{align}\label{eq:G0}
  G_m(t,x,\xi;x_0) = \xi a(x_0) (-t+ S(x;x_0))
\end{align}
with $x\to S(x;x_0)$ a global diffeomorphism from $\Rm$ to $\Rm$ with both $\partial_x S$ and $\partial_S x$ uniformly bounded and $|a(x_0)|$ bounded above and below by positive constants. We further assume that $\phi_m(x,z;x_0)$ is independent of $\xi$ as this is the case in the practical applications we consider. 

Let $(x_t,0)$ be defined as the unique point in $\Phi^{-1}(\Gamma)$ such that $t=S(x_t;x_0)$. 
Then we have the following result:
\begin{theorem}\label{thm:nondispersive}
  Assume the hypotheses of Theorem \ref{thm:localerror} as well as \eqref{eq:G0} with $\phi_m$ independent of $\xi$. Then \eqref{eq:erroru0} holds for the leading term in the wavepacket expansion $u_0$ defined as 
  \begin{align*}
     u_0(t,x,y) & = v_{\sharp0}(t,\Phi^{-1}(x,y)), \quad v_0(t,x,y) = \eps^{-\frac14}\Psi_0(t,x,z), \\ \Psi_0(t,x,z) &= \eps^{-\frac14} {\mathfrak f}_0\Big(t,x,\frac{a(x_0)}{\sqrt\eps} \big(t-S(x;x_0)\big);x_0\Big) \phi_m(x,z;x_0), \\
      {\mathfrak f}_0(t,x,X;x_0)&= \dint_{\Rm} e^{i\xi X} f_0(t,x,\xi;x_0)d\xi.
  \end{align*}
  
\end{theorem}
\begin{proof}
  This is a direct application of Theorem \ref{thm:localerror} with $f_0(t,x,\xi;x_0)\phi_m$ the leading term in the expansion, for which estimating the integral defining \eqref{eq:Psi} directly gives  the above result.
\end{proof}
This result shows that the leading term of the wavepacket in the un-rectified coordinates is concentrated in the $\sqrt\eps-$vicinity of $\Phi(x_t,0)\in\Gamma$ and propagating along $\Gamma$ with speed $|\dot x_t|=|\partial_x S(x_t;x_0)|^{-1}$.  This result generalizes constructions obtained in \cite{bal2021edge} for the Dirac operator and extends \eqref{eq:u0Dirac} to the setting of curved domain walls. 

\subsubsection{Leading term in dispersive wavepackets.}

We now make the following assumptions on $G_m(t,x,\xi;x_0)$ modeling a dispersive branch of continuous spectrum. 
\begin{itemize}
\item[(ix)] Let $\tilde \Xi \subset \Xi \subset \Rm$, all open sets, and assume $G_m(t,x,\xi;x_0)$ (real-valued) satisfies (ii) on $\Xi$. For $x_0$ fixed, consider the equation $\partial_\xi  G_m(t,x,\xi)=0$. At $(t,\xi)\in\Rm\times\Xi$ fixed, we assume the existence of a unique solution $x_m(t,\xi)$ to that equation. We denote by $\mX_t$ the range of $\Xi\ni\xi\to x_m(t,\xi)$, which we assume to be open. For $(t,x)\in \Rm\times \mX_t$, we assume the existence of a unique solution $\xi_m(t,x)$ to the same equation.

For $(t,x)\in\Rm\times \mX_t$, we assume that $|\partial^2_\xi G_m(t,x,\xi_m(t,x))|\geq C|t| |\xi_m|^a$ for $C>0$ and $a\in\Rm$. We use the notation $G''_m=\partial^2_\xi G_m(t,x,\xi_m(t,x))$.

For $(t,x,\xi)\in\Rm\times (\Rm\backslash \mX_t)\times \tilde\Xi$, we assume that $|\partial_\xi G_m(t,x,\xi)|\geq C(d(x,\mX_t)+t\aver{\xi}^b)$ for $b\in\Rm$ and $|\partial^2_\xi G_m|(t,x,\xi)\leq C \aver{t}$.
\end{itemize}
In the applications considered in this paper, $a=-3$ and $b=a+1=-2$. Above, $d(x,y)$ is Euclidean distance. With these assumptions, we have
\begin{theorem}\label{thm:dispersive}
  Assume the hypotheses of Theorem \ref{thm:localerror} as well as  (ix) above. Assume that the wavepacket $u^J$ is based on a sufficiently smooth and rapidly decaying (in $\xi$) initial condition $\hat f_m(x,\xi;x_0)$ in \eqref{eq:initcond} with support in $\tilde\Xi$ in the $\xi-$variable. Then \eqref{eq:erroru0} holds, where, using $y=\sqrt\eps z$,
  \begin{align*}
      u_0(t,x,y) &= \Phi^{-*} v_{\sharp0}(t,x,y) ,\qquad v_0(t,x,y) = \eps^{-\frac14}\Psi_0(t,x,z), \\ 
     \Psi_0(t,x,z) &= e^{\frac i{\sqrt\eps}G_m(t,x,\xi_m)}   
      e^{-i\frac\pi4 \sgn{tG''_m}}
      \Big(\frac{2\pi}{|tG''_m|}\Big)^{\frac12} f_0(t,x,\xi_m) \phi_m(x,\xi_m,z) + O_{L^\infty} \Big(\frac{\eps^{\frac14}}{|tG''_m|^{\frac12}}\Big)
  \end{align*}
  for $x\in \mX_t$ while  $|\Psi_0(t,x,z)|\leq C \aver{t}|t|^{-2} \eps^{\frac 14}$ for $x\in \Rm\backslash \mX_t$.
\end{theorem} 
The interpretation of the result is as follows. We assume an initial condition for $v_0$ concentrated in the $\sqrt\eps-$vicinity of $(x_0,0)\in \Phi^{-1}(\Gamma)$ and hence with an amplitude of order $\eps^{-\frac12}$. For any non-vanishing time, we obtain by dispersion that $v_0(t,x,y)$ has an amplitude of order $\eps^{-\frac14}|t|^{-\frac12}$ in a $\sqrt\eps$-vicinity of $\mX_t$, where $\mX_t$ is an interval including $x_0$ that grows with time $t$. Away from that interval, we obtain that $v_0$ is comparatively negligible and of order $|t|^{-1}$. The estimate on $v_0$ may in principle be improved to an arbitrary power of $\eps$ so long as $x$ remains a positive distance away from $\mX_t$ by non-stationary phase techniques; we do not pursue this further here.

In practice, we consider two settings, one in which (viii) holds with then $\tilde\Xi=\Xi=\Rm$ or $\Xi=\{|\xi|>\delta\}$ with $\delta$ arbitrary, and one where (viii) does not hold and wavepackets display dispersion only for sufficiently large wavenumbers in $\tilde\Xi$; see section \ref{sec:Dirac} for an application to the Dirac equation, with a similar result for the  Klein-Gordon equation.

\begin{proof}
  With assumption (ix), we apply a standard stationary phase argument \cite{dimassi1999spectral} to obtain the result on $\mX_t$.   
  When $x\not\in \mX_t$ and $\xi\in\tilde\Xi$, we apply a non-stationary phase argument as follows. Define the first-order differential operator
  \begin{align*}
     \rM^* = \frac{\sqrt\eps}{\partial_\xi G_m} D_\xi,\quad \rM^* e^{\frac i{\sqrt\eps} G_m} =  e^{\frac i{\sqrt\eps} G_m} ,\quad \rM = D_\xi \frac{\sqrt\eps}{\partial_\xi G_m} ,\quad \mbox{ so that:}
  \end{align*}
  \begin{align*}
     \dint_{\Rm} e^{\frac i{\sqrt\eps} G_m} f(\xi) d\xi =  \dint_{\Rm} e^{\frac i{\sqrt\eps} G_m} \rM f(\xi) d\xi 
      =  \sqrt\eps \dint_{\Rm} e^{\frac i{\sqrt\eps} G_m}\Big( -i\frac{\partial^2_\xi G_m}{(\partial_\xi G_m)^2} f - \frac{1}{\partial_\xi G_m} D_\xi f \Big) d\xi.
  \end{align*}
  Applying such an estimate with $f$ replaced by $f_0\phi_m (\xi)$, which is sufficiently smooth by assumption, provides the above estimate on $\Psi_0(t,x,z)$ for $x\not\in \mX_t$. Indeed, since $|\partial^2_\xi G_m|$ is bounded above uniformly by $\aver{t}$ and $|\partial_\xi G_m|$ is bounded below by $d(x,\mX_t) + t \aver{\xi}^b\geq t \aver{\xi}^b$ by assumption in (ix), we obtain that $\Psi_0$ in \eqref{eq:Psi0} is bounded by $\eps^{\frac14} \aver{t}|t|^{-2}$. 

 Applying the operator $\rM^k$ shows that the integral is in fact of order $\eps^{\frac k2}$ for large values of $k$ (depending on priori regularity assumptions on $f_0$) uniformly in $x$ provided that $d(x,\mX_t)$ is bounded below by a positive constant, derivatives of $G_m$ satisfy appropriate bounds (satisfied for Dirac and Klein-Gordon models), and $f$ is sufficiently smooth. We do not pursue this further here.
\end{proof}

As a corollary of the above theorem, we obtain the same dispersion result for $u_0=\Phi^{-*} v_0$ when $\Gamma$ is open. The solution $u_0(t,x,y)$ is then primarily concentrated in the $\sqrt\eps-$vicinity of the curve $\Phi(\mX_t,0)\subset\Gamma$, with an amplitude of order $\eps^{-\frac14}|t|^{-\frac12}$.

When $\Gamma$ is closed, we still have that $v_0$ disperses as written in the above result. However, in the original unrectified coordinates, $u^J$ is now well approximated by a function
\begin{align*}
   u_0(t,x,y) =  \Phi^{-*} \Big(\dsum_{k\in\Zm} v_0(t,x+k\rL,y) \Big)(t,x,y).
\end{align*}
For sufficiently small times such that $\mX_t \cap (\mX_t+\rL)=\emptyset$, then $u_0(t,x,y)$ has an amplitude bounded by $\eps^{-\frac14} |t|^{-\frac12}$. However, for longer times, where $u_0$ above is seen as a sum of a finite number of terms, we observe that  $u_0$ is still concentrated in the $\sqrt\eps$ vicinity of  (the whole of) $\Gamma$ with the amplitude still of order $\eps^{-\frac14}$. However, its amplitude no longer decays as $|t|^{-\frac12}$ and may in fact intermittently increase as $\rL-$translates of $v_0$ possibly interfere constructively.

%
%
\begin{remark}(Turning points and confinement) \label{rem:turning} \rm
Assumption (ix) above in particular requires the phase function $G_m(t,x,\xi;x_0)$ to remain real-valued. In the presence of turning points (caustics in one space dimension), $G_m(t,x,\xi;x_0)$ may develop a positive imaginary part for $x$ outside of $[x_L,x_R]\ni x_0$ with at least one of the points $x_L(\xi,x_0)$ or $x_R(\xi,x_0)$ bounded. 
Outside of $[x_L,x_R]$, the ansatz \eqref{eq:Psi} with our construction of $\psi(t,x,\xi,z)$ performed for a fixed $\xi\in\Rm$ fails to be valid. Assume $x_L$ bounded. Since $\Im G_m(t,x,\xi;x_0)>0$ for $x<x_L$, this region is forbidden classically and the wavepacket is expected to be exponentially small. As it approaches $x_L$, the wavepacket is expected to `turn around'  into a wavepacket with wavenumber $-\xi$. Modeling this coupling requires a more precise ansatz generalizing known constructions and estimate analyses for the one-dimensional Schr\"odinger equation; see \cite[Chapter 15]{hall2013quantum}, \cite{paul1993construction}, and references cited there. Alternatively, one may construct a general parametrix for solutions of $Lu=0$ as in \cite{drouot2022semiclassical}; see also \cite[Chapter 2]{guillemin1990geometric} for a corresponding analysis of the effects of turning points for the Schr\"odinger equation.

When both $x_R$ and $x_L$ are finite, we thus expect a wavepacket generated at $x_0$ with wavenumber $\xi$ to be remained confined in the interval $[x_R,x_L]$ with limited energy tunneling out of that interval. The dispersive estimate of Theorem \ref{thm:dispersive} with a  wavepacket amplitude of order $t^{-\frac12}$ therefore holds at best only for short times (until the wavepacket reaches one end point of $[x_R,x_L]$). 

\end{remark}
\section{Dirac operators} \label{sec:Dirac}

We apply the theory of section \ref{sec:wae} to a general semiclassical Dirac operator given by 
\begin{align}\label{eq:Dirac}
   L_D= \eps D_t + H_D,\qquad H_D=\eps (\gamma\sigma)\cdot(D+h) + \kappa \sigma_3,
\end{align}
where $\gamma(x,y)=(\gamma_{ij})_{1\leq i,j\leq 2}(x,y)$ is real-valued invertible  $2\times2$ matrix with uniformly (in $(x,y)$) positive determinant and $D=(D_i)_{1\leq i\leq 2}$ and $h=(h_i)_{1\leq i\leq 2}$ are seen as column vectors. In coordinates, $(\gamma\sigma)\cdot(D+h)=\gamma_{ij}\sigma_j(D_i+h_i)$ with summation over repeated indices. Above,  $\sigma_{1,2,3}$ are the standard (Hermitian) Pauli matrices and $D_{1,2}$ are identified with $D_{x,y}$. All spatially varying coefficients are assumed in $C^\infty(\Rm^2)$. 

Let $\mu(x,y)$ be a smooth density uniformly bounded above and below by a positive constant and denote by $(\cdot,\cdot)_\mu$ the inner product on the $L^2_\mu(\Rm^2;\Cm^2)$ space with measure $\mu(x,y)dxdy$.  We assume ${\rm det}\,\gamma>0$ for concreteness and the symmetry condition:
\begin{align}\label{eq:condh}
   \Im\ h = - \frac12 (\mu \gamma)^{-T}  \nabla \cdot (\mu  \gamma).
\end{align}
In \eqref{eq:Dirac},  the only term not formally of order $O(\eps)$ in \eqref{eq:Dirac} is the real-valued domain wall $\kappa(x,y)\sigma_3$. The above model therefore does not include large magnetic potentials as in \cite{bal2022magnetic,drouot2022semiclassical}, to which we refer for effects of large magnetic fields. The simplest example, generalizing the operator given in the introduction section,  is $\gamma_{ij}(x,y)=\delta_{ij}$ and $\mu(x,y)=1$ in which case we require that the imaginary part $\Im\ h=0$.

\subsection{Dirac equation}
We first obtain the following:
\begin{lemma}\label{lem:SADirac}
  Let $\gamma$ be real-valued uniformly invertible, $\kappa$ real-valued and $h$ satisfying \eqref{eq:condh}. Then $H_D=L_D-\eps D_t$ is self-adjoint as an unbounded operator on $L^2_\mu(\Rm^2;\Cm^2)$. Moreover, (vi) holds with  $q=1$ and 
$
    \mL[u] =  \sqrt{(u,u)_\mu}.
$
\end{lemma}
\begin{proof}
The constraint on $h$ written in coordinates is given by:
\begin{align*}
  \Im (\gamma^Th)_j =  \Im\ \gamma_{ij} h_i = -\frac12 \mu^{-1} \partial_i (\gamma_{ij}\mu) = -\frac12 \mu^{-1}(\nabla\cdot \gamma\mu)_j ,\qquad j=1,2.
\end{align*}
We observe by integrations by parts that 
\begin{align*}
  \dint_{\Rm^2} \bar \phi_m \gamma_{ij} \sigma_j^{mn} (D_i+h_i) \psi_n \mu d^2x 
   = \dint_{\Rm^2} -D_i(\bar\phi_m\gamma_{ij}\mu) \sigma_j^{mn}\psi_n d^2x +  \dint_{\Rm^2} \bar \phi_m \gamma_{ij} \sigma_j^{mn} h_i\psi_n \mu d^2x 
   \\ =  \dint_{\Rm^2} \overline{D_i\phi_m} \gamma_{ij} \sigma_j^{mn}\psi_n \mu d^2x - \dint_{\Rm^2} \bar\phi_m D_i(\gamma_{ij}\mu) \sigma_j^{mn}\psi_n d^2x +  \dint_{\Rm^2} \bar \phi_m \gamma_{ij} \sigma_j^{mn} h_i\psi_n \mu d^2x  \\
    =  \dint_{\Rm^2} \overline{\gamma_{ij} \sigma_j^{nm}D_i\phi_m} \psi_n \mu d^2x - \dint_{\Rm^2} \bar\phi_m \mu^{-1} D_i(\mu\gamma_{ij}) \sigma_j^{mn}\psi_n \mu d^2x +  \dint_{\Rm^2} \overline{\bar h_i\phi_m} \gamma_{ij} \sigma_j^{mn} \psi_n \mu d^2x,
\end{align*}
so that 
\begin{align*}
  (\phi,  \gamma_{ij} \sigma_j (D_i + h_i)\psi)_\mu = ( (\gamma_{ij} \sigma_j(D_i+\bar h_i)  + \mu^{-1} D_i (\gamma_{ij}\mu)\sigma_j) \phi,\psi)_\mu.
\end{align*}
When \eqref{eq:condh} holds, $H_D$ is thus formally self-adjoint  for the inner product $(\cdot,\cdot)_\mu$. That $H_D$ is self-adjoint as an unbounded operator is then shown as in \cite{thaller2013dirac}. We may therefore construct the evolution operator $e^{-itH_D}$ so that  by the Duhamel principle, $L_Du= \eps g$ with square integrable initial conditions $u(t=0)$ and square integrable source term $g(t)$ uniformly in time admits a unique solution  $u(t)\in L^2_\mu(\Rm^2;\Cm^2)$ uniformly in time $t$ over a compact interval. Moreover,
\begin{align*}
   \|u(t)\|_{L^2_\mu} \leq C( \|u(0) \|_{L^2_\mu}  + \aver{t} \sup_\tau \|g(\tau)\|_{L^2_\mu}).
\end{align*}
Here and below we use $\aver{\alpha}=\sqrt{1+|\alpha|^2}$ for $\alpha$ real- or vector-valued. This verifies assumption (vii) with $\mL[u]=\sqrt{(u,u)_\mu}$ for $q=1$.
\end{proof}


\subsection{Multiscale operator.}

Define $P = \Phi^* L_D \Phi^{-*}$ the Dirac operator in rectified coordinates. From the result $\Phi^* D \Phi^{-*}= \rJ^{-T} D$ of section \ref{sec:rectification}, we deduce that 
\begin{align*}
  P &=  \eps D_t + \eps (\gamma_P \sigma)\cdot( D + h_P)   + \kappa_P \sigma_3\\[2mm]
      \gamma_P &=  \rJ^{-1} (\gamma_L \circ \Phi) ,\qquad h_P = \rJ^{T}( h_L \circ \Phi),\qquad \kappa_P =\kappa_L \circ \Phi .
\end{align*}
By assumption on $\gamma_L$ and construction of $\rJ$, we obtain that $\gamma_P$ is uniformly invertible on $\Rm\times[-2\eta,2\eta]$ with $\det \rJ^{-1}=1$ on $\Phi^{-1}(\Gamma)$. Below we drop the subscript $P$. Associated to the operator $P$ is a multiscale operator $\rPm=\trPm \chi_\eta(\sqrt\eps z)$ for each (still to be identified) branch $m$ given by
\begin{align*}
  \trPm &=  \eps D_t + \sqrt\eps \partial_t G_m  +   (\gamma(x,\sqrt\eps z)\sigma)_1 \big(\eps D_x+\sqrt\eps \partial_x G_m + \eps h_1(x,\sqrt\eps z)\big) \\[1mm]
     & \quad + (\gamma(x,\sqrt\eps z) \sigma)_2 \big(\sqrt\eps D_z + \eps h_2 (x,\sqrt\eps z)\big) +  \kappa(x,\sqrt\eps z) \sigma_3 .
\end{align*}

We then decompose $\eps^{-\frac12} \rPm= \rPmz+\eps^{\frac12} \rPmo + \eps \rPmt$ with
\begin{align} \label{eq:rPmjD}
  \rPmz & = \partial_t G_m + (\gamma(x,0)\sigma)_1  \partial_x G_m + (\gamma(x,0)\sigma)_2 D_z + z \partial_y\kappa(x,0)\sigma_3 \\[2mm]\nonumber
  \rPmo &=  D_t +  z (\partial_y\gamma(x,0)\sigma)_1 \partial_x G_m +   (\partial_y\gamma(x,0)\sigma)_2 zD_z \\ \nonumber& \quad +  (\gamma(x,0)\sigma)_1 (D_x+h_1(x,0)) + (\gamma(x,0)\sigma)_2 h_2(x,0)  + \frac12 z^2 \partial^2_y \kappa(x,0)\sigma_3.
\end{align}
The operators $\rPmz$ and $\rPmo$ are polynomial in $z\in\Rm$ as requested in \eqref{eq:P0} and (i). The operator $\rPmt$ involves two contributions: one from $(\chi_\eta(\sqrt\eps z)-1)(\rPmz+\sqrt\eps\rPmo)$ and one from the Taylor expansion about $z=0$. In more detail, 
\begin{align} \label{eq:rPmtD}
  \eps  \rPmt =  \big(\eps^{-\frac12}\trPm-\rPmz-\sqrt\eps \rPmo)\chi_\eta(\sqrt\eps z) +(\rPmz+\sqrt\eps \rPmo)  (\chi_\eta(\sqrt\eps z)-1) .
\end{align}

\subsection{Microscopic equilibrium.} 

The leading operator is given by
\begin{align}\label{eq:rPmzD}
  \rPmz = \partial_t G_m + (\gamma\sigma)_1 \partial_x G_m + (\gamma\sigma)_2 D_z + \rho z \sigma_3,
\end{align}
where we use $\gamma(x)$ for $\gamma(x,0)$ and where we define $\rho(x):=\partial_y\kappa(x,0)$, which is strictly positive thanks to \eqref{eq:condkappa} and our convention of choice of rectified coordinates that $\kappa>0$ for $y>0$.
Define also $E=-\partial_t G_m$ and $\zeta=\partial_x G_m$. Thus $\rPmz$ may be recast as
\begin{align*}
   \rPmz= (\gamma\sigma)\cdot (\zeta, D_z) + \rho z \sigma_3 -E.
\end{align*}
The analysis of $\rPmz$ and its inverse is based on first recasting it in a normal form. To this end, we define
\begin{align}\label{eq:lambda}
   \lambda(x) = \sqrt{\gamma_{21}^2 + \gamma_{22}^2}(x),\quad \hat \sigma_2 = \lambda^{-1} (\gamma_{21}\sigma_1+\gamma_{22}\sigma_2),\quad \hat\sigma_1=\lambda^{-1}(\gamma_{22}\sigma_1 - \gamma_{21}\sigma_2).
\end{align}
Denoting by $\hat \sigma_3=\sigma_3$, we verify that $(\hat\sigma_j)_{j=1,2,3}$ are Pauli matrices satisfying the commutation relations $\hat\sigma_i\hat\sigma_j+\hat\sigma_j\hat \sigma_i=2\delta_{ij}I_{2}$. With these definitions, we have
\begin{align}\label{eq:alphabeta}
  (\gamma\sigma)_2 = \lambda\hat\sigma_2,\ \  (\gamma\sigma)_1=\alpha \hat\sigma_1+\beta\lambda\hat\sigma_2,\ \  \alpha=\frac{\gamma_{11}\gamma_{22}-\gamma_{21}\gamma_{12}}{\lambda}, \ \  \beta=\frac{\gamma_{11}\gamma_{21}+\gamma_{12}\gamma_{22}}{\lambda^2}.
\end{align}
On $\Phi^{-1}(\Gamma)$, $\det \rJ^{-1}=1$ and the matrix $\gamma$ in $L_D$ has positive determinant so that $\inf_{x\in\Rm} \alpha(x)>0$. We also have $\inf_{x\in\Rm} \lambda(x)>0$ for the same reason. We thus recast 
\begin{align*}
   \rPmz= \alpha \zeta \hat \sigma_1 + \lambda(\beta \zeta + D_z)\hat\sigma_2 + \rho z \sigma_3 -E = e^{-i\beta\zeta z} (\alpha \zeta \hat \sigma_1 + \lambda D_z\hat\sigma_2 + \rho z \sigma_3 -E )e^{i\beta\zeta z} .
\end{align*}
Let $Q_0=\frac{1}{\sqrt 2}(\sigma_3+\sigma_1)$ the global permutation of matrices such that $Q_0^*(\sigma_1,\sigma_2,\sigma_3)Q_0=(\sigma_3,-\sigma_2,\sigma_1)$. Let $\tilde Q(x)=e^{-i\frac {\check\theta(x)}2\sigma_3}$ be the standard spinorial rotation with $\check \theta(x)$ chosen based on \eqref{eq:lambda} such that $\tilde Q^*(\hat\sigma_1,\hat\sigma_2,\sigma_3)\tilde Q=(\sigma_1,\sigma_2,\sigma_3)$. Define finally $Q(x,\zeta z) = Q_0 \tilde Q(x) e^{-i\beta \zeta z}$. Then we observe that 
\begin{align}\label{eq:Q}
  Q^*\rPmz Q = H -E \qquad \mbox{ and } \qquad  Q^*(\hat\sigma_1,\hat\sigma_2,\sigma_3) Q = (\sigma_3,-\sigma_2,\sigma_1).
\end{align}

Applying the unitary  (for $L^2(\Rm_z;\Cm^2)$) transformation $Q$ thus brings us to the following normal form of the operator:
\begin{align*}
    H - E = \alpha\zeta \sigma_3 - \lambda D_z \sigma_2 + \rho z \sigma_1 - E  = \begin{pmatrix} \alpha\zeta-E & \fa \\ \fa^* & -\alpha\zeta-E \end{pmatrix},\quad  \fa=\lambda\partial_z + \rho z.
\end{align*}
All coefficients $\alpha,\lambda,\rho$ are uniformly strictly positive and bounded. We next observe that
\begin{align}\label{eq:HsquaredD}
    H^2-E^2 = (\alpha\zeta)^2 - \lambda^2 \partial_z^2 + (\rho z)^2 + \rho\lambda \sigma_3 = \begin{pmatrix} (\alpha\zeta)^2 -E^2 +\fa\fa^* & 0 \\ 0 &(\alpha\zeta)^2 -E^2 +\fa^*\fa \end{pmatrix}.
\end{align}
The eigenvalues for the harmonic oscillator $\fa^*\fa$ are $2n\lambda\rho$ for $n\geq0$. The normalized eigenvectors are $\varphi_n(x,z)$ and satisfy $\fa\varphi_n=\sqrt{2n\lambda\rho} \varphi_{n-1}$.  More precisely, let $\tilde\varphi_n(z)$ the normalized eigenvectors of $(-\partial_z+z)(\partial_z+z)\tilde\varphi_n(z)=2n\tilde\varphi_n(z)$. These are the Hermite functions given by
\begin{align*}
  \tilde\varphi_0(z)= \pi^{-\frac14} e^{-\frac12 z^2},\quad \tilde\varphi_n(z)=\frac{1}{2^{\frac n2}{\sqrt{n!}}} (-\partial_z+z)^n \tilde\varphi_0(z),\ n\geq1.
\end{align*}
We then verify that $\varphi_n$ and $\tilde\varphi_n$ are related by the following scaling change:
\begin{align}\label{eq:varphin}
  \varphi_n(x,z) = (\lambda^{-1}\rho)^{\frac14}(x) \tilde\varphi_n(\sqrt{\lambda^{-1}\rho}(x) z).
\end{align}

We thus find that the kernel of $H^2-E^2$ in  \eqref{eq:HsquaredD} is non-trivial provided that
\begin{align*}
   E_m(x,\zeta) = \epm \sqrt{(\alpha(x)\zeta)^2 + 2n \lambda(x)\rho(x)},\qquad E_0(x,\zeta)=-\alpha(x) \zeta.
\end{align*}
We consider settings where $E_m$ is real-valued whereas then,
\begin{align*}
  \alpha(x) \zeta = \pm  (E_m^2 - 2n(\lambda\rho)(x))^{\frac 12}
\end{align*}
will be real valued when $E_m^2\geq 2n(\lambda\rho)(x)$ and purely imaginary (with positive imaginary part) when $E_m^2<2n(\lambda\rho)(x)$. For $\xi\in\Xi_m$, we will obtain that $\zeta$ is indeed real-valued and non-vanishing.  The choice of sign for $\zeta$ will be given below.

We denote by $M$ the collection of indices
\begin{align}\label{eq:M}
   M = \big\{ m=(n,\epm), \ \mbox{ with } \ \epm\in\{\pm1\} \ \mbox{ for } \ n\geq1 \ \mbox{ and } \ \epm=-1 \mbox{ for } \ n=0\big\}.
\end{align}
The kernel of $H-E_m$ is then found to be
\begin{align}\label{eq:tphim}
   \tilde\phi_m(x,\zeta,z) 
   =  c_m \begin{pmatrix} \sqrt{2n\lambda\rho}  \varphi_{n-1}(x,z) \\ (E_m-\alpha\zeta) \varphi_n (x,z)\end{pmatrix},\quad  c_m = \frac{1}{\sqrt{ 2n\lambda\rho + |E_m-\alpha\zeta|^2}}.
\end{align}

Since $E=-\partial_t G_m$ and $\zeta=\partial_x G_m$, the eikonal equation for $G_m$ is then
\begin{align}\label{eq:eikonalD}
  \partial_t G_m + \epm \sqrt{(\alpha\partial_x G_m)^2 + 2n \lambda\rho } =0.
\end{align}

We consider the following global solutions to the above equation:
\begin{align}\nonumber
  G_m(t,x,\xi;x_0) &=  -\epm B_n(\xi;x_0) t + A_n(x,\xi;x_0), \quad B_n(\xi,x_0) = \sqrt{(\alpha(x_0)\xi)^2 + 2n (\lambda\rho)(x_0)} \\
   A_n(x,\xi;x_0) &= \sgn{\xi}  \int_{x_0}^x \big((\alpha(x_0)\xi)^2+ 2n ((\lambda\rho)(y)-(\lambda\rho)(x_0))\big)^{\frac12}\frac{dy}{\alpha(y)}. \label{eq:GmD}
\end{align} 
$A_n$ is constructed such that $\partial_x A_n(x_0,\xi;x_0)=\xi$. With this, we obtain that
\begin{align}\label{eq:zetaD}
  \zeta(x,\xi;x_0) = \partial_x A_n(x,\xi;x_0) = \sgn{\xi} \frac{\big((\alpha(x_0)\xi)^2+ 2n ((\lambda\rho)(x)-(\lambda\rho)(x_0))\big)^{\frac12}}{\alpha(x)}.
\end{align}
The $(\cdot)^{\frac12}$ above is the usual principal square root with branch cut on the negative real numbers and such that $(-1)^{\frac12}=i$. While $A_n(x,\xi;x_0)$ remains real-valued for $|x-x_0|$ sufficiently small, it becomes complex-valued as soon as there are points $y$ on the interval between $x_0$ and $x$ such that $(\alpha(x_0)\xi)^2+ 2n ((\lambda\rho)(y)-(\lambda\rho)(x_0))<0$. Our choice of $\xi\in\Xi_m$ in \eqref{eq:Xim} below guarantees that $\zeta$ remains real-valued. Note that when $n=0$, then $A_0$ is proportional to $\xi$ and defined for all $\xi\in\Rm$. 

The corresponding eigenvalue
\begin{align}\label{eq:EmD}
  E_m(x,\zeta(x,\xi;x_0)) = -\partial_t G_m = \epm B_n(\xi;x_0) = \epm  \sqrt{(\alpha(x_0)\xi)^2 + 2n (\lambda\rho)(x_0)}
\end{align}
is then by construction independent of $x$. The eigenvector $\tilde \phi_m(x,\zeta(x,\xi;x_0),z)$, however, depends non-trivially on $(x,\xi)$. For $m=(0,-1)$, we obtain the following simplifications
\begin{align*}
 \zeta_0(x,\xi;x_0) = \xi \frac{\alpha(x_0)}{\alpha(x)},\qquad E_0(\xi;x_0)=-\alpha(x_0)\xi, \qquad \tilde \phi_0(x,z) =  \begin{pmatrix} 0 \\ 1 \end{pmatrix}\varphi_0(x,z).
\end{align*}
This is the relativistic mode considered in \cite{bal2022magnetic,bal2021edge}.

It remains to define the following functions in the kernel of $\rPmz$:
\begin{align}\label{eq:phimD}
  \phi_m(x,\xi,z;x_0) =
  Q(x,\partial_x A_n(x,\xi;x_0) z) \tilde\phi_m(x,\partial_x A_n(x,\xi;x_0),z).
\end{align}


We recall that the coefficients $\rho$, $\lambda$, $\alpha$, and $\beta$ are defined in \eqref{eq:rPmzD}, \eqref{eq:lambda}, and \eqref{eq:alphabeta}. The function $\zeta(x,\xi;x_0)$ is defined in \eqref{eq:zetaD}. The eigenvectors $E_m$ are defined in \eqref{eq:EmD} while the eigenvectors $\phi_m$ are defined in \eqref{eq:Q}, \eqref{eq:varphin}, \eqref{eq:tphim}, and \eqref{eq:phimD}. The phase functions $G_m$ are defined in \eqref{eq:GmD}.
We now verify:
\begin{lemma}\label{lem:microscopicD}
  Let $M\ni m=(n,\epm)$ as in \eqref{eq:M} and define
  \begin{align*}
    \xim := \alpha^{-1}(x_0) \sup_{x\in\Rm} \sqrt{2n((\lambda\rho)(x)-(\lambda\rho)(x_0))}.
  \end{align*}  
  Let $x_0\in\Rm$ be fixed and define
  \begin{align}\label{eq:Xim}
      \Xi_m= \Rm \ \mbox{when } \ \xim=0,\quad \mbox{ and }\quad  \Xi_m=\{|\xi|\geq \xim+\delta\} \ \mbox{ when } \ \xim>0      ,
  \end{align}
  for some $\delta>0$.  Then assumptions (i-ii) hold for $m\in M$ and $\Xi=\Xi_m$. Moreover, $E_m$ and $G_m$ are homogeneous of degree $1$ in $(\xi,\sqrt n)$ while $\phi_m$ is homogeneous of degree $0$ in $(\xi,\sqrt n)$.
\end{lemma}
\begin{proof}
  (i) is clear. The solutions $G_m$ of \eqref{eq:eikonalD} defined in \eqref{eq:GmD} satisfy (ii) by inspection. The functions $\tilde \phi_m$ for $m\in M$  defined in \eqref{eq:tphim} are verified to form an orthonormal family of $L^2(\Rm_z;\Cm^2)$ since the rescaled Hermite functions $\varphi_n(x,z)$ form an orthonormal family of $L^2(\Rm_z;\Cm)$ as a well-known result.  This property is preserved by $Q(x,\zeta z)$ so that $\phi_m(x,\xi,z;x_0)$ also forms such a family at each $(x,\xi;x_0)$ fixed.   We thus obtain a family of solutions to $\rPmz\psi_0=0$ of the form $\psi_0(t,x,\xi,z;x_0)=f_0(t,x,\xi;x_0)\phi_m(x,\xi,z;x_0)$.
\end{proof}

\subsection{Transport equation.}

We now analyze the transport operator $\mT_m$ given by
\begin{align*}
  \mT_m f =    (\phi_m,\rPmo f\phi_m)_z =  (Q^*\phi_m,Q^*\rPmo Q f Q^*\phi_m)_z =  (\tilde \phi_m,Q^*\rPmo Q f\tilde \phi_m)_z.
\end{align*}
We then have the following results
\begin{lemma}\label{lem:trD}
  Under the hypotheses of Lemma \ref{lem:microscopicD}, the transport operator takes the form
  \begin{align}\label{eq:mTmD}
     \mT_m =  D_t + j_m(x,\xi)\alpha (x) (D_x + \tilde \theta(x,\xi))
  \end{align}
  for coefficients $j_m(x,\xi)$ and $\tilde\theta(x,\xi)$ that are homogeneous of degree $0$ in $(\xi,\sqrt n)$; see \eqref{eq:jmD} and \eqref{eq:mTmtemp} below for explicit expressions. The operator is symmetric in $L^2(\Rm_x;\Cm)$ with respect to the inner product with measure $\mu_P(x,0)dx$ where $\mu_P =  \det \rJ\, \mu_L\circ \Phi$. 
  
The solution to the equation $\mT_mf=0$ with initial conditions $f(0,\cdot)$ admits the explicit expression in \eqref{eq:invmTmD} below. More generally, the solution of $\mT_mf=g$ on $[-\rT,\rT]\times\Rm$ with initial conditions $f(0,\cdot)=\hat f(\cdot)$  satisfies the stability estimate in $\mCS_p$ given by
\begin{align}\label{eq:stabmTmD}
   \|f\phi_m\|_p \leq C_p ( \|\hat f\phi_m\|_p + \aver{\rT} \|g \phi_m\|_p),
\end{align}
for all $p\geq0$. In particular, (iv) holds.
\end{lemma}
\begin{proof}
\noindent{\bf Form of the operator.} From the above expression \eqref{eq:rPmjD} for $\rPmo$ and the fact that $Q^* (D_z Q)=-\beta\zeta$ for $Q$ defined in \eqref{eq:Q}, we first verify that 
 \begin{align*}
Q^*\rPmo Q &=  D_t + (z\partial_x G_m) (a\cdot \sigma) + z(D_z-\beta \zeta) (b\cdot \sigma) + \alpha (D_x+h_1) \sigma_3 + \alpha Q^*(D_x Q) \sigma_3\\ &\quad - \lambda\sigma_2 h_2 + \frac12 z^2 \partial^2_y\kappa \sigma_1.
\end{align*}
Here $a(x)$ and $b(x)$ are vector-valued functions such that $a_1=0$ and $b_1=0$. Indeed, $(\partial_y\gamma\sigma)_j$ for $j=1,2$ decomposes over $(\sigma_1,\sigma_2)$ so that after conjugation by $Q(x)$ (which maps $\sigma_3$ to $\sigma_1$), no contribution proportional to $\sigma_1$ is present in these terms. 
We next find from the explicit expression in \eqref{eq:tphim} that 
\begin{align*}
  (\tilde\phi_m,\sigma_1\tilde\phi_m)_{\Cm^2} & =c_m^2 2 \sqrt{2n\lambda\rho}(E_m-\alpha \Re \zeta) \varphi_{n-1}\varphi ,\qquad (\tilde\phi_m,\sigma_2\tilde\phi_m)_{\Cm^2} =0,\\
  (\tilde\phi_m,\sigma_3\tilde\phi_m)_{\Cm^2}  & =c_m^2 (2n\lambda\rho \varphi_{n-1}^2 - (E_m-\alpha\zeta)^2 \varphi_n^2).
\end{align*}
We observe from \eqref{eq:Q} that $Q^* D_xQ=Q_0^*\tilde Q^*D_x \tilde Q  Q_0 = -\frac12\check\theta'(x)Q_0^* \sigma_3 Q_0$ is proportional to $\sigma_1$ so that after product by $\sigma_3$, this is a term that decomposes over $\sigma_2$.
Therefore, using the above and $(\tilde \phi_m,z\sigma_3\tilde\phi_m)_z=(\tilde \phi_m,z^2\sigma_1\tilde\phi_m)_z=0$ by parity arguments, we find that
\begin{align*}
  (\tilde\phi_m, \big((z\partial_x G_m) (a\cdot\sigma) -\beta \zeta z (b\cdot\sigma)+ \alpha Q^*(D_x Q)\sigma_3 - \lambda\sigma_2 h_2 + \frac12 z^2 \partial^2_y\kappa \sigma_1\big) \tilde \phi_m)_z =0.
\end{align*}
As a consequence,
\begin{align*} 
  \mT_m f &= D_t f + b_3  (\tilde \phi_m, (zD_z)\sigma_3 \tilde \phi_m)_z f + (\tilde \phi_m,\alpha(D_x+h_1) \sigma_3 \tilde \phi_m)_z f .
\end{align*}
Note that $D_x$ applies to both $\tilde\phi_m$ and $f$. Let us define
\begin{align}\label{eq:jmD}
   j_m(x,\xi) = (\tilde \phi_m,\sigma_3 \tilde \phi_m)_z =  \frac{2n\lambda\rho - |E_m-\alpha\zeta|^2}{2n\lambda\rho +|E_m-\alpha\zeta|^2} = \frac{\Re (\alpha\zeta)}{E_m(\xi;x_0)} = \frac{\Re (\alpha\zeta)}{\epm B_n(\xi;x_0)}.
\end{align} 
We observe that $j_m\alpha\not=0$ for $\xi\in\Xi_m$.
The derivation of the expression for $j_m$ goes as follows. When $\zeta$ is real-valued, then 
\begin{align*}
  |E_m-\alpha\zeta|^2 &=E_m^2-2\alpha E_m\zeta + \alpha^2 \zeta^2 = 2n\lambda\rho + 2\alpha\zeta(\alpha\zeta-E_m) \\
   &= - 2n\lambda\rho +2E_m(E_m-\alpha\zeta),
\end{align*}
so that $j_m=\alpha\zeta/E_m$ since $(E_m-\alpha\zeta)\not=0$. For completeness, we note that when $\zeta$ is purely imaginary, then
$
  |E_m-\alpha\zeta|^2 = 2n\lambda\rho +\alpha^2|\zeta|^2 =\alpha^2\zeta^2+2n\lambda\rho  +\alpha^2|\zeta|^2=2n\lambda\rho
$
so that $j_m=0$ then. The hypothesis $\xi\in\Xi_m$ in \eqref{eq:Xim} implies that $\zeta$ is real-valued and non-vanishing.

There is no reason for the term $(\tilde \phi_m,zD_z \sigma_3 \tilde \phi_m)_z$ to vanish a priori and in fact it should not  in the presence of genuine anisotropy so that $\mT_m$ preserves its symmetric structure.  We find more specifically that 
\begin{align*}
   (\tilde\phi_m,z\partial_z\sigma_3\tilde\phi_m)_z &=c_m^2 \dint_{\Rm}[ 2n\lambda\rho \varphi_{n-1} z\partial_z \varphi_{n-1} - (E_m-\alpha\zeta)^2 \varphi_{n} z\partial_z \varphi_{n} ]dz\\
    &= c_m^2 \dint_{\Rm}[ 2n\lambda\rho \tilde\varphi_{n-1}(z) z\partial_z \tilde\varphi_{n-1}(z) - (E_m-\alpha\zeta)^2 \tilde\varphi_{n} z\partial_z \tilde\varphi_{n} ]dz\\
     &= -\frac12 c_m^2 (2n\lambda\rho - (E_m-\alpha\zeta)^2) = -\frac12 j_m(x,\xi)
\end{align*}
using for the first line: \eqref{eq:varphin}, the change of variables $\sqrt{\lambda^{-1}\rho}z\to z$, and that $z\partial_z$ is homogeneous of degree $0$ in $z$; and for the second line the well-known identities:
\begin{align*}
   z\tilde\varphi_n = \sqrt{\frac{n+1}2} \tilde\varphi_{n+1} + \sqrt{\frac{n}2} \tilde\varphi_{n-1},\qquad \partial_z\tilde\varphi_n = -\sqrt{\frac{n+1}2} \tilde\varphi_{n+1} + \sqrt{\frac{n}2} \tilde\varphi_{n-1}
\end{align*}
to deduce that the integral of $\tilde\varphi_{n} z\partial_z \tilde\varphi_{n}$ over $z$ equals $-\frac12$ independent of $n\geq1$. This allows us to recast the transport operator as 
\begin{align}\label{eq:mTmtemp}
  \mT_m &= D_t  + b_3\frac i2 j_m  + \alpha (\tilde \phi_m,\sigma_3 (D_x\tilde \phi_m))_z  + \alpha j_m (D_x+h_1)  .
\end{align}
This may be rewritten as \eqref{eq:mTmD} since for $\xi\in\Xi_m$, we observe from the definition of $\zeta(x,\xi)$ that $|j_m|>0$. Moreover, $|j_m|$ is uniformly (in $(x,\xi)$) bounded below by a $\delta-$dependent positive constant  when $|\xi|>\delta>0$ and uniformly bounded below by a constant times $|\xi|$ when $\Xi_m=\Rm$.

\medskip
\noindent{\bf Symmetry of the operator.}
The operator $\mT_m$ is symmetric in the following sense.
 We have $P=\Phi^{*} L_D \Phi^{-*}$ so that 
\begin{align*}
  (\psi_1, L \psi_2 )_{\mu_L} = (\psi_1, \Phi^{-*} P \Phi^{*} \psi_2)_{\mu_L} = (\Phi^*\psi_1,P \Phi^*\psi_2)_{\mu_P}
\end{align*}
with $\mu_P=\det \rJ\  \mu_L\circ \Phi$. Therefore,
\begin{align*}
  (\Phi^*\psi_1,P \Phi^*\psi_2)_{\mu_P} = (\psi_1, L_D \psi_2 )_\mu = (L_D \psi_1, \psi_2 )_\mu =   (P \Phi^*\psi_1, \Phi^*\psi_2)_{\mu_P}
\end{align*}
so that $P$ is symmetric with inner product $\mu_P$. This holds for functions $\psi_j$ compactly supported in any domain where $\Phi$ is a bijection.

We drop the subscript $P$ in $\mu$ with now $\mu=\mu_P$.  Since $G_m$ is real-valued, the terms proportional to $\partial_t G_m$ and $\partial_x G_m$ are real-valued local multiplications and therefore symmetric for the $L^2_\mu(\Rm^2;\Cm^2)$ inner product with weight $\mu_P(x,\eps^{-\frac12}z)$. Since $\phi_m$ is subgaussian in the $z$ variable, we thus obtain from the above that up to possibly a term smaller than $C\eps^N$ for each $N$, then (recall that $\rPm=\trPm \chi_\eta(\sqrt\eps z)$ with $\trPm$ a differential operator)
\begin{align*}
  (f_1\phi_m, \rPm f_2\phi_m)_\mu  =  (\rPm f_1\phi_m,  f_2\phi_m)_\mu.
\end{align*}
Note that ${\rm det}\rJ=1+\varkappa(x)\sqrt\eps z$ and thus is equal to $1$ to leading order in $\eps$. The leading term of $\rPm$ vanishes when applied to $\phi_m$ by construction. This shows that $\rPmo$ is symmetric and
\begin{align*}
  (f_1,\mT_m f_2)_x = (f_1\phi_m, \rPmo f_2\phi_m)  =  (\rPmo f_1\phi_m,  f_2\phi_m) = (\mT_m f_1,f_2)_x
\end{align*}
with $(\cdot,\cdot)_x$ the inner product in $L^2_\mu(\Rm_x;\Cm)$ with weight $\mu_P(x,0)$.  

Let us introduce $\epsilon=\epm \sgn{\xi}$ so that $j_m=\epsilon |j_m|$ as one readily verifies. 
Recall that $|j_m|>0$ is uniformly bounded below when $|\xi|>\delta>0$ and uniformly bounded below by a constant times $|\xi|$ when $\Xi_m=\Rm$. We then recast the transport operator as
\begin{align*}
  \mT_m=D_t + \epsilon \tal  (D_x +\tilde\theta) ,\qquad  \tal = |j_m|\alpha,
\end{align*}
where $\tilde\theta$ is implicit from the expression \eqref{eq:mTmtemp}.
The above symmetry implies as in the derivation of \eqref{eq:condh} that 
\begin{align*}
  \mT_m = D_t + \epsilon \tal (D_x -\frac i2 (\tal\mu)^{-1} \partial_x (\tal \mu)) + \epsilon  \tal \theta = D_t + \epsilon \tal(\tal\mu)^{-\frac12}D_x (\tal\mu)^{\frac12}  + \epsilon   \tal \theta
\end{align*}
with $\theta$ implicitly defined by the above equation and now real-valued.

\medskip
\noindent{\bf Inversion of the operator.}
Consider the solution of $\mT_m f=0$ with initial condition $f(0,x)$. Then define 
\begin{align*}
  g(t,x) =  e^{i\int_0^x\theta(y)dy} (\tal\mu)^{\frac12}  f (t,x), \quad \mbox{ so that } \quad  (D_t +  \epsilon\tal  D_x  )g=0,
\end{align*}
as one verifies. 
Introduce next the change of variable $dy/dx=\tal^{-1}(x)$ which, thanks to the above properties of $|j_m|>0$, generates a diffeomorphism $\ry(x)$ from $\Rm$ to $\Rm$ preserving orientation and with inverse $\rx(y)$.  We observe that 
\begin{align*}
  (D_t+\epsilon D_y) (g\circ x) (t,y)=0,\quad  (g\circ x)(t,y)= (g\circ x)(0,y-\epsilon t).
\end{align*}
Thus 
\begin{align*}
   g(t,x) = g(0,\rx(\ry(x)- \epsilon t))
\end{align*}
and the solution to $\mT_m f(t,x)=0$ with initial condition $f(0,x)$ is given explicitly by
\begin{align}\label{eq:invmTmD}
   f(t,x) = (\tal\mu)^{-\frac12}(x) e^{-i\int_0^x \theta(y)dy} \big[ (\tal\mu)^{\frac12}(\cdot) e^{i\int_0^{\cdot} \theta(y)dy}  f(0,\cdot) \big] (\rx(\ry(x)-\epsilon t)).
\end{align}
We thus obtain an explicit expression for the solution of the equation $\mT_m f=0$ with initial condition $f(0,x)$ when $\xi\in\Xi_m$. Applying a standard Duhamel principle provides a solution to the equation $\mT_mf=g$ with vanishing initial conditions as well.

\medskip
\noindent{\bf Stability of the inverse operator.}  The coefficients involved in the explicit expression \eqref{eq:invmTmD} are all homogeneous of degree $0$ in $(\xi,\sqrt n)$ thanks to Lemma \ref{lem:microscopicD} and hence so are their derivatives with respect to $x$. All bounds therefore are independent of $|\xi|\geq\delta$. We recall that $\theta(x,\xi)$ is real-valued so that the complex exponentials and their derivatives with respect to the variable $x$ are all bounded uniformly in $x$. This shows that derivatives of order $p$ in $(t,x)$ of $f(t,x)$ are controlled by derivatives of $f(0,x)$ of order up to $p$. We thus obtain that 
\begin{align*}
  \|f(t,x,\xi)\|_p \leq C_p \|f(0,x,\xi)\|_p
\end{align*}
with $C_p$ independent of the dependence of the above functions in $\xi$ when $\Xi_m \cap \{|\xi|<\delta\}=\emptyset$. 

The case $\Xi_m=\Rm$ may be considered when $\lambda\rho$ is constant. The transport equation when $n=0$ poses no difficulty since $|j_m|$ is then bounded from below uniformly in $\xi$. Consider $n\geq1$ and $\lambda\rho$ constant. In this case, $\alpha(x)\zeta(x,\xi)=\alpha(x_0)\xi$ in \eqref{eq:zetaD} is independent of $x$ and the dependence in $x$ of $\tilde \phi_m$ is solely in the components $\varphi_n(x,z)$ in \eqref{eq:varphin}. Since $j_m$ in \eqref{eq:jmD} is then independent of $x$ as well, then $\partial_x j_m=0$ implies that $(\tilde\phi_m,\sigma_3 \partial_x \tilde\phi_m)_z=0$ since $\sigma_3$ is hermitian. This implies that the transport operator is given by
\begin{align*}   
  \mT_m &= D_t  + j_m(\xi)( b_3(x)\frac i2  + \alpha(x)(D_x+h_1(x)) )  = D_t + j_m(\xi) \alpha(x)(D_x+\tilde\theta(x))
  \\ \nonumber
  &= D_t + j_m(\xi) \alpha^{\frac12}\mu^{-\frac12} D_x (\alpha\mu)^{\frac12} + j_m(\xi) \alpha(x) \theta(x) 
\end{align*}
with $\theta(x)$ real-valued and independent of $\xi$. For $\xi=0$, we thus have $j_m=0$ and hence $\mT_m=D_t$ whose inversion is stable. Introducing the change of variables $dy/dx=\alpha^{-1}(x)$ now independent of $\xi$, the solution of the equation $\mT_mf=0$ is given by
\begin{align}\label{eq:invmTmD1}
   f(t,x) = (\alpha\mu)^{-\frac12}(x) e^{-i\int_0^x \theta(y)dy} \big[ (\alpha\mu)^{\frac12}(\cdot) e^{i\int_0^{\cdot} \theta(y)dy}  f(0,\cdot) \big] (\rx(\ry(x)- j_m(\xi) t)).
\end{align}
This explicit expression shows that derivatives of order $p$ in $(t,x)$ of $f(t,x)$ are controlled by derivatives of $f(0,x)$ of order up to $p$ independent of $j_m(\xi)$ for $|\xi|\leq\delta$. 

By the Duhamel principle, we observe that for solutions of the equation $\mT_m f=g$ with initial condition $f(0,x,\xi)$, we also obtain in either case of $\Xi_m$ defined in \eqref{eq:Xim} that 
\begin{align*}
  \|f(t,x,\xi)\|_p \leq C_p \big( \|f(0,x,\xi)\|_p + \aver{t} \|g(t,x,\xi)\|_p).
\end{align*}
In practice, it is convenient to use \eqref{eq:stabmTmD}, which, as one readily verifies, is a direct consequence of the above estimate using the explicit expression for the norms $\|\cdot\|_p$. This proves condition (iv).
\end{proof}


\begin{remark}\label{rem:transportD}\rm
   The transport operator $\mT_m$ was shown to be stably invertible for both forms of $\Xi_m$ in \eqref{eq:Xim}. Wavepackets may be chosen arbitrarily as described in Theorem \ref{thm:globalerror} when $\lambda\rho$ is constant.
   
    When $\lambda\rho$ is not constant, then $\alpha\zeta$ in \eqref{eq:zetaD} becomes purely imaginary for values of $(x,\xi)$ such $\alpha^2(x_0)\xi^2+2n(\lambda\rho)(x_0)-2n(\lambda\rho)(x)<0$. If $x_0$ is not chosen as a maximum of $\lambda\rho$, then such an inequality holds for $|\xi|$ small enough. At such points $j_m(x,\xi)=0$ in \eqref{eq:jmD}. Furthermore, at a fixed value of $\xi$, $j_m$ is real valued on an interval $(x_-,x_+)\ni x_0$ with $x_-$ or $x_+$ finite. At such a point, say $x_-$, we observe when $|\partial_x(\lambda\rho)(x_-)|>0$ that $j_m$ approaches $0$ with a scaling of the form  $\sqrt{x-x_-}$. The change of variables $dy/dx=\tal^{-1}$ no longer generates a diffeomorphism from $\Rm$ to $\Rm$, an expression of the form \eqref{eq:invmTmD} ceases to be valid for arbitrary times, and we no longer have any guarantee of stability of the transport solution. 
   
   Points $(x,\xi)$ with $\zeta(x,\xi;x_0)$ purely imaginary correspond to phase functions $G_m$ with positive imaginary parts and hence tunneling regions that are forbidden classically; see our discussion in remark \ref{rem:turning}. As wavepackets approach turning points such as $x_-$ or $x_+$ above, we expect wavepackets to turn around and remain in the energetically favorable domain $(x_-,x_+)$. The ansatz \eqref{eq:Psi} based on the analysis of the propagation of a single mode $m$ with global phase $G_m$ in \eqref{eq:GmD} ceases to be valid and the transport solution as constructed above should not be controlled.
\end{remark}

We now move to the analysis of the operator $\rPmzi$. 
\begin{lemma}
   Under the hypotheses of Lemma \ref{lem:microscopicD}, let $N_m=N_m(x,\xi;x_0)$ be the one-dimensional space of spinors in $L^2(\Rm_z;\Cm^2)$ proportional to $Q(x,\zeta z) (\fa\varphi_n,(E_m-\alpha\zeta)\varphi_n)^T$ with $\zeta=\zeta(x,\xi;x_0)$. Let $p=(p_1,p_2)$. For $g\in \mCS_{(p_1,p_2+1)} \cap N_m^\perp$, the equation $\rPmz f=g$ admits solutions in $\mCS_p$ given by $f=\rPmzi g + f_1\phi_m$ for an arbitrary scalar-valued function $f_1(t,x,\xi;x_0)$ such that $f_1\phi_m\in \mCS_p$. The solution operator $\rPmzi g$ provides the unique solution in $\mCS_{p} \cap N_m^\perp$. Explicitly, we have
\begin{align}\label{eq:rPmziD}
   \rPmzi &= Q (H-E_m)^{-1} Q^*,\\ \nonumber
   (H-E_m)^{-1} &= \begin{pmatrix} (\fa\fa^*-2n\lambda\rho)^{-1} &0 \\ 0 & (\fa^*\fa-2n\lambda\rho)^{-1} \end{pmatrix}\begin{pmatrix} \alpha\zeta+E_m & \fa \\ \fa^* & E_m-\alpha\zeta \end{pmatrix}.
\end{align}
  In particular, assumption (iii) holds with $q_2=1$.
\end{lemma}
\begin{proof} We analyze the equation
\begin{align*}
    \rPmz \psi = g ,\quad\mbox{ or equivalently} \quad Q^* \rPmz Q Q^*\psi = (\partial_t G_m+H)  Q^*\psi = Q^* g.
\end{align*}
Recall that $E_m=-\partial_t G_m$ and look for equation $(H-E_m)\psi=g$ instead (using $g$ and $\psi$ instead of $Q^*g$ and $Q^*\psi$ to simplify), which implies $(H^2-E_m^2)\psi=(H+E_m)g=:f$.

From \eqref{eq:HsquaredD}, we observe that $(H^2-E_m^2)\psi=f$ admits a solution if and only if $(\fa\varphi_n,f_1)=(\varphi_n,f_2)=0$.
To solve $(H-E_m)\psi=g$, we have $f=(H+E_m)g$ so that 
\begin{align*}
(\fa\varphi_n,f_1)&=(\fa\varphi_n,(\alpha\zeta+E_m)g_1+\fa g_2)=(\alpha\zeta+E_m)(\varphi_n,\fa^* g_1) + (E_m^2-(\alpha\zeta)^2)(\varphi_n,g_2)\\ &= (\alpha\zeta+E_m) [(\varphi_n,\fa^*g_1) + (E_m-\alpha\zeta) (\varphi_n,g_2) ]= (\alpha\zeta+E_m) (\varphi_n,f_2).
\end{align*}
Since $\alpha\zeta+E_m\not=0$ unless $n=0$, we deduce that $(\fa\varphi_n,f_1)=0$ when $(\varphi_n,f_2)=0$ for $f=(H+E_m)g$. Therefore, $(H-E_m)\psi=g$ admits solutions if and only if $g$ satisfies the $m-$dependent condition
\begin{align*} 
   (\varphi_n, \fa^* g_1 + (E_m-\alpha\zeta) g_2)= (\fa\varphi_n,g_1) + (E_m-\alpha\zeta) (\varphi_n,g_2) =0.
\end{align*}
Let $N_m=N_m(x,\zeta)$ be the space of spinors $\psi$ that are proportional to $(\fa\varphi_n,(E_m-\alpha\zeta)\varphi_n)^T$ so that the above condition means that $g$ is orthogonal to $N_m$ for the usual inner product in $L^2(\Rm_z;\Cm^2)$. On $N_m^\perp$, we define $(H-E_m)^{-1}$ as in \eqref{eq:rPmziD}.
This inverse is defined also when $\zeta\in i\Rm$ so long as $E_m$ is real-valued.  In terms of regularity, the inverse is stable in $(z,D_z)$ and multiplies by $\aver{\zeta}$ or equivalently $\aver{\xi}$ in that variable. All these functions are smooth in $x$. The operators $(\fa\fa^*-2n\lambda\rho)^{-1}$ and $(\fa^*\fa-2n\lambda\rho)^{-1}$ are bounded on their respective domains of definition since the different simple eigenvalues of $\fa^*\fa$ are separated by multiples of $2\lambda\rho$, which is uniformly bounded below by a positive constant. On the other hand, the remaining operators are homogeneous of degree $1$ in $(\xi,\sqrt n)$ for the diagonal terms and of the form $\fa$ or $\fa^*$ for the off-diagonal terms. Each of these operators is bounded from $\mCS_{(p_1,p_2+1)}$ to $\mCS_{p}$ by construction of the norms. This property remains after conjugation by $Q(x)$. This proves \eqref{eq:boundrPmzi} and concludes the proof of the lemma.
\end{proof}
\begin{lemma}
  The operators $\rPmj$ are bounded from $\mCS_{(p_1+j,p_2+1)}$ to $\mCS_p$  and in particular, (v) holds.
\end{lemma}
\begin{proof}
  We look at the explicit expressions in \eqref{eq:rPmjD} and \eqref{eq:rPmtD}. Since $D_x$, $D_t$, and multiplication by $\xi$ or $\zeta$ are bounded from $(p_1,p_2+1)$ to $p$ and $D_z$ and $z$ are bounded from $(p_1+1,p_2)$ to $p$, we directly obtain the result for $\rPmz$ and $\rPmo$. We do not need to assume that $\xi\in\Xi_m$ here since $\zeta(x,\xi)$ is defined for all $\xi\in\Rm$. For the bound on $\rPmt$, we need to analyze two terms. One comes from the Taylor expansion in $\sqrt\eps z$ given by
  \begin{align*}
    f(\sqrt\eps z) = f(0) + \sqrt\eps z f'(0) +  (\sqrt\eps z)^2 \int_0^1 (1-t) f''(\sqrt\eps zt) dt.
  \end{align*}
  We thus observe a multiplication by $z^j$ in $\rPmj$ resulting in a bound from  $(p_1+2,p_2)$ to $p$ for that first contribution in \eqref{eq:rPmtD}. The second contribution may be recast as
  \begin{align*}
     \frac{\chi_\eta (\sqrt\eps z)-1}{\eps \aver{z}^2}\aver{z}^2\rPmz+ \frac{\chi_\eta (\sqrt\eps z)-1}{\sqrt\eps \aver{z}} \aver{z}\rPmo. 
  \end{align*}
  Since $\chi_\eta(y)=0$ in the vicinity of $y=0$, we observe that $(\eps^j\aver{z}^j)^{-1}(\chi_\eta (\sqrt\eps z)-1)$ is uniformly bounded with uniformly bounded derivatives. We then obtain that both $\aver{z}^2\rPmz$ and $\aver{z}\rPmo$ are bounded in the prescribed sense. 
\end{proof}
%
%
%
%
%

\subsection{Wavepacket analysis.} 
\label{sec:diracwa}
Since hypotheses (i)-(vi) hold with $\Xi=\Xi_m$, we may apply Theorem \ref{thm:localerror} to obtain solutions $u$ of $Lu=0$ with arbitrary accuracy in powers of $\eps$ in the sense of $L^2_\mu(\Rm^2;\Cm^2)$. 

\paragraph{Relativistic mode.}
When $m=(n=0,\epm=-1)$, we observe that 
\begin{align*}
   G_{(0,-1)} (t,x,\xi;x_0) = \alpha(x_0)\xi  \Big(t + \int_{x_0}^x \frac{dy}{\alpha(y)} \Big).
\end{align*}
Thus for this mode, (i)-(v) always hold with $\Xi=\Rm$ and  we may apply Theorem \ref{thm:nondispersive} to obtain that wavepackets are concentrated in the vicinity of $x_t$ defined by
\begin{align*}
  t+\int_{x_0}^{x_t} \frac{dy}{\alpha(y)} =0.
\end{align*}
\begin{theorem}\label{thm:propzero}
  When $m=(0,-1)$, we obtain the following expressions for the leading term of the propagating wavepackets
  \begin{align*}
   \Psi_0(t,x,z) &= \eps^{-\frac14}  {\mathfrak f}_0\Big(t,x,\frac{\alpha(x_0)}{\sqrt\eps} \big(t + \dint_{x_0}^x  \alpha^{-1}(y) dy \big)\Big) \varphi_0(z)  Q(x) \begin{pmatrix} 0\\1\end{pmatrix},\\
   v_0(t,x,y) &= \eps^{-\frac12}  {\mathfrak f}_0\Big(t,x,\frac{\alpha(x_0)}{\sqrt\eps} \big(t + \dint_{x_0}^x  \alpha^{-1}(z) dz\big)\Big) \varphi_0\Big(\frac{y}{\sqrt\eps}\Big)  Q(x) \begin{pmatrix} 0\\1\end{pmatrix},\\
    u_0(t,x,y) &=v_0(t,\Phi^{-1}(x,y)),
  \end{align*}
where we have defined ${\mathfrak f}_0(t,x,X)=\int_{\Rm} e^{iX\xi} f_0 (t,x,\xi)d\xi$ for $\hat f_0(x,\xi;x_0)$ arbitrary (sufficiently smooth). 
\end{theorem}
Recall that $\hat f_0$ is the prescribed initial condition \eqref{eq:initcond}  in the construction of the wavepacket $\psi^J_0$. The above result is written with $\Gamma$ open with an obvious modification when $\Gamma$ is closed since $f_0(t,x,\xi;x_0)$ has support in $x$ concentrated in the vicinity of $x_t$.
\paragraph{Arbitrary initial conditions.\!\!} 
We now prove the following result on a construction of a parametrix.
\begin{lemma} \label{lem:arbitraryDirac}
 Assume:
 \begin{enumerate}
   \item[(vii)$_D$]  The function $\lambda\rho$ is constant, with $\rho(x)=\partial_y\kappa(x,0)$ and $\lambda(x)=\sqrt{\gamma_{21}^2(x)+\gamma_{22}^2(x)}$ for $x\in\Rm$.
\end{enumerate}
Then (vii)-(viii) hold with $\Xi=K=\Rm$ independent of $m$. In particular, Theorem \ref{thm:globalerror} applies.
\end{lemma}
(vii)$_D$ holds for instance when $\gamma$ is constant and the strength of the domain wall $|\nabla\kappa|$ remains constant along $\Gamma$. For an arbitrary tensor $\gamma(x)$ (smooth and with uniformly positive determinant) and hence an arbitrary coefficient $\lambda(x)$, we can always construct a domain wall $\kappa(x,y)$ (with fixed $\Gamma=\kappa^{-1}(0)$) such that (vii)$_D$ holds.
\begin{proof}
(vii)$_D$ implies (vii) and (viii) hold for each $m$ with $\Xi=K=\Rm$. Indeed, 
\begin{align}\label{eq:changeAn}
  A_n (x,\xi;x_0) = \xi \int_{x_0}^x  \frac{\alpha(x_0)}{\alpha(y)}dy = (x-x_0)k =: \xi(x-x_0)\Lambda(x,x_0)
\end{align}
with $\Lambda(x,x_0)$ bounded above and below by strictly positive constants independent of $m$. The completeness hypothesis in (vii) with $q=1$ comes from the completeness of the Hermite functions in $L^2(\Rm)$ and the explicit expression \eqref{eq:tphim} while (viii) is then clear since $q=1$. Thus, Theorem \ref{thm:globalerror} applies: any (smooth) initial condition appropriately localized in the vicinity of a point $x_0\in\Gamma$ may be decomposed over wavepackets whose propagation is described at the macroscopic level.
\end{proof}

\paragraph{High frequency initial conditions. } 

A similar change of variables to the one in \eqref{eq:changeAn} also applies for sufficiently high-frequency initial wavepackets. How high-frequency is, however, $m-$dependent. For $x\in \Xi_m$, the change of variables 
\begin{align*}
  k_n(x,\xi) = \frac{A_n(x,\xi;x_0)}{x-x_0}
\end{align*}
defines a diffeomorphism from $\Rm\times\Xi_m$ to its image for each fixed $x_0$ as one easily verifies. We can use this change of variables to show that any initial condition in mode $m$ with sufficiently high-frequency content $k_n$ may be represented as a wavepacket of the form \eqref{eq:Psi}. We leave the details to the reader. 

\paragraph{Dispersive modes. } 

Let us consider a first result that applies when (vii)$_{D}$ holds.
\begin{lemma}
  When (vii)$_{D}$ holds, then (ix) holds as well with $\tilde\Xi=\Xi=\Rm$ and Theorem \ref{thm:dispersive} applies.
\end{lemma}
\begin{proof}
  Let $\alpha_0=\alpha(x_0)$. Then, $G_m(t,x,\xi;x_0)=-\epm \sqrt{\alpha_0^2\xi^2 + 2n \lambda\rho} + \xi \int_{x_0}^x \frac {\alpha_0}{\alpha(y)}dy$ so that 
  \begin{align*}
     \partial_\xi G_m = -\epm t \frac{\alpha_0^2\xi}{\sqrt{\alpha_0^2\xi^2 + 2n \lambda\rho}} + \int_{x_0}^x \frac {\alpha_0}{\alpha(y)}dy,\qquad 
     \partial^2_\xi G_m = \frac{-\epm t 2n (\lambda\rho) \alpha_0^2}{(\alpha_0^2\xi^2 + 2n \lambda\rho)^{\frac32}}.
  \end{align*}
  Since $\alpha$ is uniformly bounded above and below by positive constants, there is a unique solution $x_m(t,\xi)$ of $\partial_\xi G_m(t,x,\xi;x_0)=0$. Define the points $x_{t\pm}$ via
  \begin{align*}
    \dint_{x_0}^{x_{t\pm}} \alpha^{-1}(y) dy = \pm |t| . 
  \end{align*}
  We observe that $x_m(t,\xi)\in \mX_t=(x_{t-},x_{t+})$ for all $\xi\in\Rm$ and that $x_{t\pm}$ are asymptotically attained as $|\xi|\to\pm\infty$. For $(t,x)$ fixed, we deduce from the constant sign of $\partial^2_\xi G_m$ that there is at most one solution $\xi_m(t,x)$ such that $\partial_\xi G_m(t,x,\xi;x_0)=0$. So for $(t,x)\in \Rm\times\mX_t$, we have a unique solution $\xi_m$ to that equation. The bound on $\partial^2_\xi G_m$ in (ix) for  $(t,x)\in \Rm\times\mX_t$ is clear from its explicit expression with $a=-3$.
  
  Finally, assume that $(t,x)\in\Rm\times (\Rm\backslash \mX_t)$. Then, from the expression of $\partial_\xi G_m$ and the fact that $\alpha^{-1}$ is bounded above and below by positive constants, we obtain that 
\begin{align*}
  |\partial_\xi G_m(t,x,\xi)|  &=  |\partial_\xi G_m(t,x,\xi) - \partial_\xi G_m(t,x_m,\xi) | \geq C |x-x_m| \\
  & \geq C \inf_\pm |x-x_{t\pm}|+|x_{t\pm}-x_m|  \geq  C d(x,\mX_t) + C |t| \aver{\xi}^{-2},
\end{align*}
since, using the above explicit expression for $\partial_\xi G_m(t,x_m,\xi)=0$,
\begin{align*}
 \inf_\pm |x_{t\pm}-x_m| &\geq \inf_\pm \Big| \int_{x_m}^{x_{t\pm}} \alpha^{-1}(y)dy \Big|  = \inf_\pm  \Big| \Big( \int_{x_0}^{x_{t\pm}} -\int_{x_0}^{x_m} \Big)\alpha^{-1}(y)dy \Big|  
   \\ & =  |t|\Big( 1- \frac{\alpha_0|\xi|}{\sqrt{\alpha_0^2\xi^2+2n\lambda\rho}}\Big) \geq C |t| \aver{\xi}^{-2}.
\end{align*}
This concludes the derivation of (ix) with $b=-2$ so that Theorem \ref{thm:dispersive} applies.
\end{proof}
\begin{remark}\rm
  As in Theorem \ref{thm:dispersive}, the above dispersive estimates apply to $u_0=\Phi^{-*} v_{0\sharp}$ with $v_0(t,x,y)=\Psi_0(t,x,z)$. We obtain a decay in $|t|^{-\frac12}$ when $\Gamma$ is open and in both cases $\Gamma$ open or closed, a reduction of the amplitude of $u_0$ by $\eps^{\frac14}$ between time $0$ and arbitrary positive times for all points $x\in \mX_t$. Outside of $\mX_t$, the solution $u_0$ is smaller with a maximal amplitude reduced by (at least) another order $\eps^{\frac14}$ for non-zero times.
  
  We observe that each wavepacket within a narrow frequency band concentrated about a given $\xi$ propagates at speed given by 
\begin{align} \label{eq:sppedm}
   \dot x_m(t,x) = \pdr{}t x_m(t,x) = \dfrac{\epm \alpha_0 \alpha(x_m) \xi}{\sqrt{\alpha_0^2\xi^2+2n \lambda\rho}},
\end{align}
for $x\in \mX_t$.
When $n=0$, we retrieve that $\dot x=-\alpha(x)$ independent of $\xi$ for the relativistic mode.  For $n\geq1$, we also observe that modes with large wavenumber $\xi$ are close to relativistic with a speed that converges to $\epm\alpha(x)$ as $|\xi|\to\infty$.

Repeated iterations of the non-stationary phase argument presented in the proof of Theorem \ref{thm:dispersive} show that the decay of $\Psi_0$ outside of $\mX_t$ is faster than any power of $\eps$. 
\end{remark} 

\paragraph{Local dispersion effect.} The preceding result assumes $\Xi=\Rm$ in (ix), which only holds for the specific model with $\lambda\rho$ constant. When $\lambda$ and $\rho$ are arbitrary, we saw that the component $A_n(x,\xi;x_0)$ may have a non-trivial imaginary component. However, when $\xi\in\Xi_m$ defined in \eqref{eq:Xim}, no such effect occurs and dispersion may be estimated. We consider this setting here.
\begin{enumerate}
   \item[(ix')$_D$] Assume $n$ fixed and recall that 
  $
     \xim = \alpha_0^{-1} \sup_{x\in\Rm} \sqrt{2n((\lambda\rho)(x)-(\lambda\rho)(x_0))}.
   $
    Define $\Xi_\delta = \{\xi\in\Rm; |\xi| > \xim+\delta\}$ for $\delta>0$. 
    We assume (ix) with  $\Xi=\Xi_\delta$ and $\tilde\Xi=\Xi_{2\delta}\subset\Xi_\delta$.
\end{enumerate}
\begin{lemma}\label{lem:localdispersion}
  Under (ix')$_D$,  the results stated in Theorem \ref{thm:dispersive} apply.
\end{lemma}
\begin{proof}
   We assume that $\xi\in\Xi\ = \Xi_\delta$ implying that $A_n(x,\xi;x_0)$ is real-valued. Define $\beta(x)=2n(\lambda\rho)(x)$ and $\beta_0=\beta(x_0)$ while $\tilde\beta(x)=\beta(x_0)-\beta(x)$. With $x_0$ fixed, we have
   \begin{align*}
      G_m(t,x,\xi) = -\epm t \sqrt{\alpha_0^2\xi^2+2n(\lambda\rho)(x_0)} + \sgn{\xi}  \dint_{x_0}^x \sqrt{\alpha_0^2\xi^2  + \tilde\beta(y)} \frac1{\alpha(y)}dy.
   \end{align*}
   We use repeatedly for $a^2>0$ and $b\geq0$ that
   \begin{align*}
      \partial_\xi \sqrt{a^2\xi^2+b} = \frac{a^2\xi}{\sqrt{a^2\xi^2+b}},\qquad \partial_\xi \sqrt{a^2\xi^2+b} = \frac{a^2b}{(a^2\xi^2+b)^{\frac32}}.
   \end{align*}
   Then,
   \begin{align*}
    \partial_\xi G_m = \frac{-\epm t\alpha_0^2\xi}{\sqrt{\alpha_0^2\xi^2+2n(\lambda\rho)(x_0)}} + \frac{\xi^2}{|\xi|}  \dint_{x_0}^x \frac{\alpha_0^2}{\sqrt{\alpha_0^2\xi^2  + \tilde\beta(y)} }\frac1{\alpha(y)}dy.
   \end{align*}
   For $|\xi|>\xim+\delta$, the above integrand is strictly positive and bounded so that there is a unique $x_m(t,\xi)$ such that $\partial_\xi G_m(t,x_m(t,\xi),\xi)=0$. Moreover, at that point $x=x_m(t,\xi)$, we find after some algebra that 
   \begin{align*}
     |t| = \Big| \dint_{x_0}^x \sqrt{1+\frac{\beta(y)}{\alpha_0^2\xi^2+\tilde\beta(y)}} \frac1{\alpha(y)}dy \Big| \leq C |x-x_0| \aver{ (|\xi|-\xim)^{-1}}
   \end{align*}
   so that the signal propagates with a group velocity at least equal to $C\delta$ for some constant $C$.  
   
   An additional differentiation yields
   \begin{align*}
      \partial^2_\xi G_m = \frac{-\epm t\alpha_0^2\beta(x_0)}{(\alpha_0^2\xi^2+\beta(x_0))^{\frac32}} + \sgn{\xi}  \dint_{x_0}^x \frac{\alpha_0^2\tilde\beta(x_0)}{(\alpha_0^2\xi^2  + \tilde\beta(y))^{\frac32} }\frac1{\alpha(y)}dy.
   \end{align*}
   At the point $x$ where $\partial_\xi G_m=0$, this gives
   \begin{align*}
      \sgn{\xi} \partial^2_\xi G_m &=  \frac{-\beta(x_0)}{\alpha_0^2\xi^2+\beta(x_0)}  \dint_{x_0}^x \frac{\alpha_0^2}{(\alpha_0^2\xi^2  + \tilde\beta(y))^{\frac12} }\frac1{\alpha(y)}dy
+  \dint_{x_0}^x \frac{\alpha_0^2\tilde\beta(x_0)}{(\alpha_0^2\xi^2  + \tilde\beta(y))^{\frac32} }\frac1{\alpha(y)}dy\\
 &= \dint_{x_0}^x \frac{\-\alpha_0^4 \xi^2 \tilde\beta(y) }{(\alpha_0^2\xi^2+\beta(x_0))(\alpha_0^2\xi^2  + \tilde\beta(y))^{\frac32} }\frac1{\alpha(y)}dy
   \end{align*}
   from which we deduce that at such points, 
   \begin{align*}
      |\partial^2_\xi G_m (t,x,\xi)| \geq C \aver{\xi}^{-3} |x-x_0| \geq C \aver{\xi}^{-3} |t| \aver{(|\xi|-\xim)^{-1}}.
   \end{align*}
   This proves the lower bound for $G''_m$ for $x\in \mX_t$ and $\xi\in\Xi=\Xi_\delta$.
   
   It remains to obtain a lower-bound for $\partial_\xi G_m$ when $x\not\in \mX_t$ and $\xi\in \tilde\Xi=\Xi_{2\delta}$. Assume $\xi>0$ for concreteness. For $(t,\xi)\in\Rm\times\Xi_{2\delta}$, there is a unique $x_m(t,\xi)$ such that $\partial_\xi G_m(t,x_m,\xi)=0$. Recall that $\mX_t=(x_-,x_+)$ is defined as the range of $x_m(t,\xi)$ for $\xi\in\Xi=\Xi_\delta$ and is an interval for $\xi\geq\xim+\delta$. Define $\tilde\mX_t=(\tilde x_-,x_+)$  as the range of $x_m(t,\xi)$ for $\xi\in\tilde\Xi=\Xi_{2\delta}$. We have $x_-<\tilde x_-<x_+$ for $t\not=0$.
   
   For $x\geq x_+$, we find as in the preceding case that 
   \begin{align*}
      |\partial_\xi G_m(t,x,\xi)| \geq C d(x,x_+) + C|t| \aver{\xi}^{-2}.
   \end{align*}
   For $x<x_-$, we obtain from the above expression for $\partial_\xi G_m$, the definition of $\xim$, and the fact that $\partial_\xi G_m(t,x_m,\xi)=0$, that
   \begin{align*}
      |\partial_\xi G_m(t,x,\xi)| 
        = \Big| \dint_{x}^{x_m}  \!\!\frac{\xi\alpha_0^2}{\sqrt{\alpha_0^2\xi^2  + \tilde\beta(y)} }\frac{dy}{\alpha(y)} \Big|
   \geq C\delta |x-x_m| = C \delta (|x-x_-| + |x_--\tilde x_-|).
   \end{align*}
   It remains to find a lower bound on $|x_--\tilde x_-|$ to obtain a lower bound on $|\partial_\xi G_m|$ for $\xi\in \Xi_{2\delta}$.  At fixed $t$, let $x_\tau$ be the solution of 
   \begin{align*}
      t=\int_{x_0}^{x_\tau}  \frac{\sqrt{\tau^2+\beta(x_0)}}{\sqrt{\tau^2+\tilde\beta(y)}} \frac{dy}{\alpha(y)} ,\quad \alpha_0 (\xim+\delta)\leq \tau \leq \alpha_0(\xim+2\delta).
   \end{align*}
   Differentiating with respect to $\tau$ yields
   \begin{align*}
      0= \pdr{x_\tau}\tau   \frac{\sqrt{\tau^2+\beta(x_0)}}{\sqrt{\tau^2+\tilde\beta(x_\tau)}} \frac{1}{\alpha(x_\tau)} +  \int_{x_0}^{x_\tau}  \frac{-\beta(y)}{(\tau^2+\beta(x_0))^{\frac12}(\tau^2+\tilde\beta(y))^{\frac32}} \frac{dy}{\alpha(y)} 
   \end{align*}
   from which we observe that
   \begin{align*}
      \Big|\pdr{x_\tau}\tau \Big|    \geq C \sqrt{\tau^2+\tilde\beta(x_t)}   \int_{x_0}^{x_\tau}  dy = C \delta(x_\tau-x_0) \geq C \delta^2 |t|.
   \end{align*}
   Integrating over an interval of length $\alpha_0\delta$, we find that $|x_--\tilde x_-|\geq C\delta^3 |t|$.
    Thus overall,
    \begin{align*}
      |\partial_\xi G_m(t,x,\xi)| \geq C \delta \big( (x_--x) + \delta^3 t \big).
    \end{align*}
    It remains to apply Theorem \ref{thm:dispersive} to conclude. 
\end{proof}

%
\section{Klein-Gordon operators}\label{sec:KG}
%

This section considers the general semiclassical Klein-Gordon equation
\begin{align}\label{eq:KG}
 L_{KG}=   -c^{-2}(\eps D_t)^2 + \eps D \cdot \gamma \eps D + \kappa^2 - \eps V = \eps^2c^{-2}\partial_t^2 - \eps^2 \nabla \cdot \gamma \nabla + \kappa^2 - \eps V.
\end{align}
Here $\gamma$ is a uniformly positive definite symmetric tensor at each point $(x,y)\in\Rm^2$, $\kappa$ is the usual domain wall, $c$ is a uniformly positive sound speed, and $V$ is an additional real-valued coefficient satisfying the constraint that $-\eps^2 \nabla \cdot \gamma \nabla + \kappa^2 - \eps V$ remains a non-negative operator. All coefficients are assumed smooth and bounded.

A first objective is to find an expression for $V$ allowing for the presence of relativistic modes such as \eqref{eq:u0KG}.  In the operator $L_{KG}$ defined in \eqref{eq:unpKG} in the introduction, we chose $V=1$.  When the domain wall is curved, the choice of $V$ is more complicated.

On $\Gamma$ are defined the tangent vector $\tau$ to the curve $\Gamma$ and normal vector $\nu$ to it; see Figure \ref{fig:geom}.  These vectors are naturally extended to the vicinity of $\Gamma$ given by $(x,y)\in \Gamma_{2\eta}:= \Phi(\Rm\times [-2\eta,2\eta])$ by defining $(\tau,\nu)$ at $(x,y)=\Phi(\tx,\ty)$ to be the same vectors as those defined at $\Phi(\tx,0)$. In such a vicinity, we decompose $\gamma$ such that 
\begin{align*}
  D\cdot \gamma D = D_\tau(\gamma_{11}-\frac{\gamma_{12}^2}{\gamma_{22}} )D_\tau + (D_\nu+ D_\tau \frac{\gamma_{12}}{\gamma_{22}})\gamma_{22} (D_\nu+  \frac{\gamma_{12}}{\gamma_{22}} D_\tau ).
\end{align*}
Here, $D_\tau=\tau\cdot D$ and $D_\nu=\nu\cdot D$. 
We next define the operators
\begin{align*}
  \fa_\tau = \sqrt{\gamma_{11}-\frac{\gamma_{12}^2}{\gamma_{22}}}\, \sqrt \eps \partial_\tau, \quad \fa_\nu = \sqrt{\gamma_{22}}  \, \sqrt\eps(\partial_\nu+  \frac{\gamma_{12}}{\gamma_{22}} \partial_\tau )  + \eps^{-\frac12}\kappa \chi_\eta.
\end{align*}
We still denote by $\chi_\eta\in C^\infty_c(\Rm^2)$ a function equal to $1$ in the $\eta-$vicinity of $\Gamma$ and equal to $0$ for points that are a distance at least $2\eta$ away from $\Gamma$. We find that 
\begin{align*}
   L_{KG} = -c^{-2}(\eps D_t)^2 +\eps \fa_\tau^* \fa_\tau + \eps\fa_\nu^* \fa_\nu +(1-\chi_\eta^2) \kappa^2 + \sqrt{\gamma_{22}}  \, \eps(\partial_\nu+  \frac{\gamma_{12}}{\gamma_{22}} \partial_\tau )  (\kappa \chi_\eta) - \eps V.
\end{align*}
We then choose 
\begin{align}\label{eq:V}
   V = (\partial_\nu+  \frac{\gamma_{12}}{\gamma_{22}} \partial_\tau )  \kappa \chi_\eta  - W
\end{align}
and assume that $W\geq0$. An operator $L_{KG}$ with specific weight $W=0$ in the vicinity of $\Gamma$ allows for relativistic modes propagating along $\Gamma$. As soon as $W>0$, however, these modes acquire a non-vanishing mass term and disperse. The term $V$ takes the above form in the vicinity of $\Gamma$; its expression is irrelevant away from $\Gamma$ and this is the reason why it was multiplied by $\chi_\eta$.

In a variational formulation, we use the above decomposition close to $\Gamma$ and the original expression for $L$ away from $\Gamma$ to construct a meaningful energy functional. More precisely, define the energy functional (for $u$ real-valued)
\begin{align}\label{eq:energyKG}
   {\mathcal E}[u](t) = \mL^2[u](t) = \frac12 \dint_{\Rm^2} \eps c^{-2} (\partial_t u)^2 + \eps \gamma \nabla u \cdot \nabla u + \eps^{-1}\kappa^2 u^2 - V u^2 \ dxdy. 
\end{align}
We may decompose this integral as performed over $\Gamma_{2\eta}$ and $\Rm^2\backslash \Gamma_{2\eta}$. The second term has the above expression with $V=-W\leq0$ as expected from an energy functional while the integral over $\Gamma_{2\eta}$ takes the form
\begin{align*}
\frac12 \dint_{\Gamma_{2\eta}} \eps c^{-2} (\partial_t u)^2 + \eps (\fa_\tau u)^2 + \eps (\fa_\nu u)^2 + \eps^{-1}(1-\chi_\eta^2)\kappa^2 u^2 + \eps W u^2 \ dxdy. 
\end{align*}
We observe by standard integrations by parts that for $u$ sufficiently smooth and rapidly decaying:
\begin{align*}
   \dint_{\Rm^2} \partial_t u \eps^{-1} L_{KG}u \, dxdy = \dr{}t {\mathcal E}[u](t).
\end{align*}
This shows the energy conservation for any problem of the form $Lu=g$ with appropriate prescribed  initial conditions $u(t=0)$ and $\sqrt\eps \partial_t u(t=0)$. By construction, $-V\geq0$ where $\chi_\eta=0$. On the domain where $\chi_\eta=1$, i.e., in the vicinity of $\Gamma$, the integrand in ${\mathcal E}$ is given by 
\begin{align*}
   \frac12 \dint_{\chi_\eta=1}  c^{-2} (\sqrt\eps\partial_t u)^2 +  (\fa_\tau u)^2 +  ( \fa_\nu u)^2 +  W u^2 \ dxdy. 
\end{align*}

The above provides a uniqueness and existence result for the Klein-Gordon wave equation. For a source problem $L_{KG}u=g$,  we obtain
\begin{align*}
    \partial_t u \eps^{-1} Lu = \partial_t u \eps^{-1} g  = \eps^{-\frac32} g \sqrt\eps \partial_t u.
\end{align*}
We thus deduce a control of $\mL[u]$ by $\eps^{-\frac32}g$ with $-\frac32=\frac12-q$. We thus obtained:
\begin{lemma}
 Let $W\geq0$ in \eqref{eq:V}. Then (vi) holds for $\mL$ defined in \eqref{eq:energyKG} with $q=2$.
\end{lemma}

%
\paragraph{Multiscale operator.}

Consider the pulled back operator $P=\Phi^* L \Phi^{-*}$ given by
\begin{align*}
 P=   -c_P^{-2}(\eps D_t)^2 + \eps D \cdot \gamma_P \eps D + \kappa_P^2 - \eps V_P ,
\end{align*}
where the pulled-back coefficients are found to be:
\begin{align*}
  c_P=c_L\circ \Phi, \quad  \kappa_P=\kappa_L\circ \Phi, \quad  V_P=V_L\circ \Phi, \quad \gamma_P=\rJ^{-1}\gamma_L\circ \Phi \rJ^{-T}.
\end{align*}
We drop the subscript $P$ below.
 This gives the multiscale operator
\begin{align*}
  \trPm &= -c^{-2}(\eps D_t+\sqrt\eps \partial_t G_m)^2 +(\eps D_x+\sqrt \eps\partial_x G_m,\sqrt\eps D_z) \gamma  (\eps D_x+\sqrt \eps\partial_x G_m,\sqrt\eps D_z)^T \\ & \quad + (\kappa^2 - \eps V )(x,\sqrt\eps z)
\end{align*}
where $\gamma$ is also seen as a function of $(x,\sqrt\eps z)$. From this, we deduce
\begin{align*}
   \rPmz &= -c^{-2}(x) (\partial_t G_m)^2 +( \partial_x G_m,D_z) \gamma(x,0) (\partial_x G_m,D_z)^T + (\rho(x)z)^2 - V_0(x).
\end{align*}
Here, as in the Dirac case, $\rho(x)=\partial_n\kappa(x,0)$. We also write $c(x)=c(x,0)$. Finally,
\begin{align*}
   V_0(x) = \sqrt{\gamma_{22}} (x) \rho(x)
\end{align*}
is the leading term in $V=V_0+\sqrt\eps V_1$ since $\partial_\tau\kappa=0$ along $\Gamma$.

Define $E=-\partial_t G_m$ and $\zeta=\partial_x G_m$. We recast the above operator as 
\begin{align*}
   &-c^{-2} E^2 + \zeta^2 \gamma_{11} + 2\gamma_{12}\zeta D_z + \gamma_{22} D_z^2 + (\rho z)^2-V_0 \\ 
    = & -c^{-2} E^2 + \zeta^2 (\gamma_{11}- \frac{\gamma_{12}^2}{\gamma_{22}}) + \gamma_{22} (D_z + \frac{\gamma_{12}}{\gamma_{22}} \zeta)^2 + (\rho z)^2-V_0 \\
    = & \ e^{-i\frac {\gamma_{12}}{\gamma_{22}}  \zeta z} \Big(  -c^{-2} E^2 + \zeta^2 (\gamma_{11}- \frac{\gamma_{12}^2}{\gamma_{22}}) + \gamma_{22} D_z^2 + (\rho z)^2-V_0 \Big) e^{i\frac {\gamma_{12}}{\gamma_{22}} \zeta z} \\
     = & \ e^{-i\frac {\gamma_{12}}{\gamma_{22}}  \zeta z} \Big(  -c^{-2} E^2 + \zeta^2 (\gamma_{11}- \frac{\gamma_{12}^2}{\gamma_{22}}) + \fa^*\fa  + W \Big) e^{i\frac {\gamma_{12}}{\gamma_{22}} \zeta z},\qquad \fa=\sqrt{\gamma_{22}} \partial_z + \rho z.
\end{align*}
Define similarly to the Dirac case $\lambda=\sqrt{\gamma_{22}}$ and $\alpha=\sqrt{\gamma_{11}- \frac{\gamma_{12}^2}{\gamma_{22}}}$.

The eigenvalues of $\fa^*\fa$ are $2n\lambda\rho$ for $n\geq0$ and are associated with the eigenvector $\varphi_n(x,z)$ defined in \eqref{eq:varphin}. The solutions of $\rPmz\tilde \phi_m=0$ (in the variable $\zeta$) are thus given by
\begin{align*}
  E_m(x,\zeta) = c(x) \epm \sqrt{\alpha^2 \zeta^2  + 2n \lambda\rho + W},\quad \tilde \phi_m(x,\zeta,z)=e^{-i\frac {\gamma_{12}}{\gamma_{22}} \zeta z} \varphi_n(x,z).
\end{align*}
The associated Eikonal equations are thus given by
\begin{align}\label{eq:eikonalKG}
   \partial_t G_m + c(x) \epm \sqrt{\alpha^2 (\partial_x G_m)^2  + 2n \lambda\rho + W} =0.
\end{align}
These are essentially the same eikonal equations \eqref{eq:eikonalD} as for the Dirac operator. We assume $W=0$ to simplify and to allow for the presence of relativistic modes. When $W$ does not vanish, we simply replace $2n\lambda\rho$ below by $2n\lambda\rho +W$. Explicitly, we find
\begin{align}\nonumber
  G_m(t,x,\xi)&=-\epm B_n(\xi)t + A_n(x,\xi),\quad B_n(\xi) = c(x_0)\sqrt{\alpha^2(x_0)\xi^2 + 2n(\lambda\rho)(x_0)} \\
  A_n(x,\xi) &= \sgn{\xi} \dint_{x_0}^x \sqrt{\frac{c^2(x_0)}{c^2(y)} \big(\alpha^2(x_0)\xi^2+2n(\lambda\rho)(x_0)\big) -2n(\lambda\rho)(y)} \frac{dy}{\alpha(y)} \label{eq:GmKG}
\end{align}
where we assume that $\xi\in\Xi_m$ with 
\begin{align}\label{eq:XimKG}
   \Xi_m = \Big\{ \xi\in\Rm, \ \  \inf_{x\in\Rm} \frac{c^2(x_0)}{c^2(x)} \big(\alpha^2(x_0)\xi^2+2n(\lambda\rho)(x_0)\big) -2n(\lambda\rho)(x) > \delta^2  \Big\}.
\end{align}
Up to a universal constant, the definition of $\Xi_m$ in \eqref{eq:XimKG} is similar to that of \eqref{eq:Xim} in the Dirac case. Unlike the Dirac case, we do not consider the case $\delta=0$ even when $c^2\lambda\rho$ is constant since, as we will observe, the inversion of the transport operator is unstable when $n=0$ for $|\xi|$ small.


The modes $n=0$ for $\eps_m=\pm1$ are relativistic with different directions of propagation, unlike the Dirac case where only $\eps_m=-1$ is allowed.
The general solution of $\rPmz\psi_0=0$ is 
\begin{align*}
    \psi_0(t,x,\xi,z;x_0) = f_0(t,x,\xi;x_0)  \phi_m(x,\xi,z;x_0),\quad \phi_m(x,\xi,z;x_0):= e^{-i\frac {\gamma_{12}}{\gamma_{22}} z \partial_x A_n(x,\xi)} \varphi_n(x,z).
\end{align*}

Since the Hermite functions $\varphi_n$ form an orthonormal basis of $L^2(\Rm_z;\Cm)$ and $\epm=\pm1$ for $m=(n,\epm)$, we obtain that the families $\phi_m$ for $m\in M_1=\{\epm=1\}$ and $m\in M_2=\{\epm=-1\}$ also form orthonormal basis of $L^2(\Rm_z;\Cm)$. This proves (i-ii) as well as (vii) provided $\Xi_m$ may be chosen independent of $m$ as in section \ref{sec:Dirac}.  Note that when the anisotropy coefficient $\gamma_{12}=0$, then $\phi_m$ does not depend on the multiscale parameter $\xi$ nor on the location of the initial wavepacket $x_0$.

Using the above expression for $E_m$, the equation
$
   \rPmz \psi = g
$ 
is recast as
\begin{align*}
  (\fa^*\fa - 2n\lambda\rho) e^{i \frac{\gamma_{12}}{\gamma_{22}}z \partial_x A_n(x,\xi)} \psi = e^{i \frac{\gamma_{12}}{\gamma_{22}}z \partial_x A_n(x,\xi)} g .
\end{align*}

We define the kernel
\begin{align*}
  N_{x} = \Cm \phi_m(x, \partial_x A_n(x,\xi),z) = \Cm \varphi_n(x,z)  \subset L^2(\Rm_z;\Cm)
\end{align*}
and $N^\perp$ the orthogonal complement in $L^2(\Rm_z;\Cm)$. For $g\in N^\perp\cap \mCS_{p+1}$ for each $(x,\xi;x_0)$, there is a unique solution in $\mCS_p$ to $\rPmz \psi=g$ denoted by 
\begin{align} \label{eq:rPmziKG}
  \rPmzi g = e^{ -i \frac{\gamma_{12}}{\gamma_{22}}z \partial_x A_n(x,\xi)}  (\fa^*\fa - 2n\lambda\rho)^{-1}  e^{i \frac{\gamma_{12}}{\gamma_{22}}z\partial_x A_n(x,\xi)} g.
  \end{align}
Any solution of the above equation is of the form $\psi=\rPmzi g + f \phi_m$ for $f(x,\xi;x_0)\in \mCS_p(\Rm;\Cm)$. 

The operator $ (\fa^*\fa - 2n\lambda\rho)^{-1}$ is bounded since successive eigenvalues of $\fa^*\fa$ are separated by $2\lambda\rho$, which is uniformly bounded below. This proves (iii) with $q_2=0$.

\paragraph{Transport equation.}
We now look at the next-order term
\begin{align*}
  \rPmo &=  - z \partial_y c^{-2} (\partial_t G_m)^2 - 2c^{-2} \partial_t G_m D_t + (\partial_x G_m,D_z) z\partial_y\gamma  (\partial_x G_m,D_z) ^T \\ & + (D_x,0)\gamma (\partial_x G_m,0) ^T+ (\partial_x G_m,0)  \gamma (D_x,0)^T +  \partial^2_y\kappa \partial_y\kappa z^3  
    - z \partial^2_{xy}\kappa - z \partial^2_y\kappa + z \partial_y W.
\end{align*}
The transport operator is defined by and verifies the following identities
\begin{align*}
 (f_1\phi_m, \rPmo f_0\phi_m) = (f_1,\mT_m f_0)_x =  (\rPmo f_1\phi_m, f_0\phi_m) =  (\mT_m f_1, f_0)_x
\end{align*}
which shows that $\mT_m$ is symmetric for the usual inner product on $L^2(\Rm_x)$. Above, $(\cdot,\cdot)$ is the usual inner product on $L^2(\Rm^2)$. 

All terms proportional to $z$ or $z^3$ cancel when integrating in the $z$ variables in the construction of $\mT_m$ since $\varphi_{n-1}^2$ is even. Moreover, we observe that 
\begin{align*}
 \int  dz\  \partial_x G_m  \bar\phi_m z \partial_y \gamma_{12} D_z \phi_m + \overline{D_z \phi_m}  z \partial_y  \gamma_{12}  \partial_x G_m   \phi_m = 0 ,\\
 \int dz\   \partial_x G_m  \bar\phi_m z \partial_y \gamma_{11} \partial_x G_m  \phi_m + \overline{D_z \phi_m}  z \partial_y  \gamma_{22}  D_z \phi_m =0
\end{align*}
since, using $\zeta=\partial_x A_n(x,\xi)$,
\begin{align*}
  D_z \phi_m = (\frac {\gamma_{12}}{\gamma_{22}}\zeta \varphi_n -D_z \varphi_n) e^{-i\frac {\gamma_{12}}{\gamma_{22}} \zeta z}, \quad D_z\varphi_n + \overline{ D_z\varphi_n}=0.
\end{align*}
As a consequence,
\begin{align}\label{eq:mTmKG}
   \mT_m & =  -2 c^{-2} \partial_t G_m D_t + (\phi_m ,[(D_x,0)\gamma (\partial_x G_m,0) ^T+ (\partial_x G_m,0)  \gamma (D_x,0)^T] \phi_m)_z \\ \nonumber
    & = -2 c^{-2} B_n(\xi) \epm \Big(  D_t +j_m (D_x + \tilde\theta)\Big),
\end{align}
where we have defined
\begin{align*}
   j_m(x,\xi) = - \epm c^2(x) \gamma_{11}(x) \frac{\partial_x A_n(x,\xi) }{B_n(\xi)}
\end{align*}
and $\tilde\theta(x,\xi)$ implicitly. 

From the construction of $G_m$, we obtain that $B_n$ and $A_n$ are homogeneous of degree one in $(\xi,\sqrt n)$, and hence so are all coefficients in $\mT_m$ since $j_m$ and $\tilde\theta$ as constructed are homogeneous of degree zero in $(\xi,\sqrt n)$ . This contrasts with the Dirac cases where such coefficients are homogeneous of degree zero in $(\xi,\sqrt n)$ and reflects the second-order nature of the Klein-Gordon operator. 

When $n\geq1$, then $c^{-2} B_n(\xi)$ is bounded below by positive constants independent of $(x,\xi)$. However, when $n=0$, $B_n(\xi)$ is proportional to $|\xi|$ and is arbitrarily small for $|\xi|$ small. It is the reason for the choice of $\Xi_m$ in \eqref{eq:XimKG} with $\delta>0$. An equation of the form $\mT_mf=g$ is not stable when $n=0$ uniformly in $\xi$ unless $|\xi|\geq C\delta$ for $C\delta>0$, which is imposed in \eqref{eq:XimKG}.

Note that as in the Dirac setting, $|j_m|$ is uniformly bounded below by a positive constant for $\xi\in\Xi_m$ when $\delta>0$.  Up to a multiplication by $2 c^{-2} B_n(\xi)$, and a slightly different definition of the current $j_m$, the operator $\mT_m$ is the same as the one obtained in \eqref{eq:mTmD} for Dirac operators. We thus obtain the same explicit inversion procedure as the one given in Lemma \ref{lem:trD}. The main difference is that $\mT_m^{-1}$ now involves a multiplication by $(2 c^{-2} B_n(\xi))^{-1}$, which  is uniformly bounded provided $|\xi|\geq C\delta>0$. 

We thus obtain that for $\Xi$ defined as $\Xi_m$ in \eqref{eq:XimKG}, then  (iv) holds. From the definition of the operators $\rPmj$, (v) holds with $q=2$. Summarizing the above results, we thus obtained:
\begin{lemma}
 Let $\Xi=\Xi_m$ in \eqref{eq:XimKG} and $W\geq0$. Then conditions (i-vi) hold with $q=2$ and $q_2=0$. 
\end{lemma}

\begin{remark}\label{rem:transportKG}\rm
   The comments in Remark \ref{rem:transportD} still apply for the Klein-Gordon model. In particular, when $c^2\lambda\rho$ is not constant, then $G_m$ may become complex valued when $\xi\not\in \Xi_m$ in \eqref{eq:XimKG}. As for the Dirac model, solutions to the transport should no longer be controlled to reflect the fact that the ansatz \eqref{eq:Psi} is no longer meaningful to describe tunneling effects; see also remark \ref{rem:turning}.
   
   The operator $\mT_m$ in \eqref{eq:mTmKG} is similar to its counterpart for the Dirac equation with the main difference that its coefficients are bounded below by $C\aver{\xi}$ for $\xi$ large. As a result $\mT_m^{-1}$ involves a regularization of the form $\aver{\xi}^{-1}$. The operator $\mT_m^{-1}\rPmo \rPmzi\rPmo$  involves multiplication by $\aver{\xi}^{2q-1}=\aver{\xi}^3$ instead of $\aver{\xi}^{2q}$. Accounting for this improved regularity would require us to separate the role of $\aver{\xi}$ and $(D_t,D_x)$ in the definition of the spaces $\mCS_p$, which we do not do here. 
\end{remark}

We may therefore apply Theorem \ref{thm:localerror} and obtain accurate solutions of the equation $L_{KG}u=0$ by means of wavepackets of the form $u^J$. Such wavepackets satisfy the results stated in Theorem \ref{thm:nondispersive} when $n=0$ and Theorem \ref{thm:dispersive} when $n\geq1$. Hypothesis (ix) is verified exactly as in section \ref{sec:Dirac} for the Dirac operator. The main difference between the Dirac and Klein-Gordon models is the presence of two counter-propagating relativistic modes for the topologically trivial Klein-Gordon model and only one relativistic mode for the topologically non-trivial and time-reversal-symmetry-breaking Dirac model. Note that the Klein-Gordon model admits relativistic modes only for the specific choice of $V$ given by $W=0$ whereas relativistic modes naturally always appear for the Dirac model.

\medskip

We conclude this section with the following result:
\begin{lemma}
 Assume
\begin{enumerate}
  \item[(vii)$_{\rm KG}$] $c^2\lambda\rho$ is constant and $\Xi=\{ |\xi|\geq\delta\}$ for $\delta>0$.
\end{enumerate}
Then (vii)-(viii) hold and Theorem \ref{thm:globalerror} applies. 
\end{lemma}
We recall that $c^2(x)$, $\lambda(x)=\sqrt{\gamma_{22}(x)}$ and $\rho(x)= \partial_y\kappa(x,0)$ are defined along the interface in the rectified coordinates for $x\in\Rm$.
\begin{proof}
That (vii) follows from (vii)$_{\rm KG}$ is proved as for the Dirac operator since 
\begin{align*}
    A_n(t,\xi) =\xi \int_{x_0}^x \frac{c(x_0)\alpha(x_0)}{c(y)\alpha(y)} dy = (x-x_0) k :=\xi(x-x_0)\Lambda(x,x_0)
\end{align*}
with $\Lambda(x,x_0)$ bounded above and below by strictly positive constants independent of $m$. The completeness condition is again a consequence of the completeness of the Hermite functions in $L^2(\Rm)$.
To prove (viii), we observe that $B_{n1}=-B_n$ while $B_{n2}=+B_n$ with $B_n$ given in \eqref{eq:GmKG}. The $2\times2$ Vandermonde matrix $A$ in (viii) therefore has determinant $2B_n$, which is homogeneous of degree $1$ in $(\xi,\sqrt n)$ and hence uniformly bounded from below when $\xi\in \Xi$ given in (vii)$_{\rm KG}$.  Thus Theorem \ref{thm:globalerror} applies and we obtain that arbitrary initial conditions with frequencies supported in $K$ (avoiding an arbitrarily small $\eps-$independent vicinity of the point $k=0$) may be represented as a superposition of wavepackets constructed in Lemma \ref{lem:psimJ}.
\end{proof}

\section*{Acknowledgment.} The author would like to acknowledge many stimulating discussions on the topic of wavepacket propagation with Simon Becker and Alexis Drouot. This research was partially supported by the U.S. National Science Foundation, Grants DMS-1908736, DMS-2306411, and EFMA-1641100.

{\small 

\begin{thebibliography}{10}


\bibitem{avron1994}
{\sc J.~E. Avron, R.~Seiler, and B.~Simon}, {\em Charge deficiency, charge
  transport and comparison of dimensions}, Comm. Math. Phys., 159 (1994),
  pp.~399--422.

\bibitem{B-BulkInterface-2018}
{\sc G.~Bal}, {\em Continuous bulk and interface description of topological
  insulators}, Journal of Mathematical Physics, 60 (2019), p.~081506.

\bibitem{B-EdgeStates-2018}
\leavevmode\vrule height 2pt depth -1.6pt width 23pt, {\em Topological
  protection of perturbed edge states}, Communications in Mathematical
  Sciences, 17 (2019), pp.~193--225.

\bibitem{bal2022topological}
\leavevmode\vrule height 2pt depth -1.6pt width 23pt, {\em Topological
  invariants for interface modes}, Communications in Partial Differential
  Equations, 47(8) (2022), pp.~1636--1679.

\bibitem{bal2023topological}
\leavevmode\vrule height 2pt depth -1.6pt width 23pt, {\em Topological charge
  conservation for continuous insulators}, Journal of Mathematical Physics, 64
  (2023), p.~031508.

\bibitem{bal2022magnetic}
{\sc G.~Bal, S.~Becker, and A.~Drouot}, {\em Magnetic slowdown of topological
  edge states}, Accepted in Comm. Pure Applied Math. and arXiv preprint
  arXiv:2201.07133,  (2023).

\bibitem{bal2021edge}
{\sc G.~Bal, S.~Becker, A.~Drouot, C.~F. Kammerer, J.~Lu, and A.~Watson}, {\em
  Edge state dynamics along curved interfaces}, Accepted in SIAM Math. Anal.
  and arXiv:2106.00729,  (2021).

\bibitem{bal2023asymmetric}
{\sc G.~Bal, J.~G. Hoskins, and Z.~Wang}, {\em Asymmetric transport
  computations in dirac models of topological insulators}, Journal of
  Computational Physics, 487 (2023), p.~112151.

\bibitem{bal2022multiscale}
{\sc G.~Bal and D.~Massatt}, {\em Multiscale invariants of {Floquet}
  topological insulators}, Multiscale Modeling \& Simulation, 20 (2022),
  pp.~493--523.

\bibitem{BES94}
{\sc J.~Bellissard, A.~van Elst, and H.~Schulz-Baldes}, {\em The noncommutative
  geometry of the quantum {H}all effect}, Journal of Mathematical Physics, 35
  (1994), pp.~5373--5451.

\bibitem{Be13}
{\sc B.~A. Bernevig}, {\em Topological Insulators and Topological
  Superconductors}, Princeton University Press, 2013.

\bibitem{delplace2017topological}
{\sc P.~Delplace, J.~Marston, and A.~Venaille}, {\em Topological origin of
  equatorial waves}, Science, 358 (2017), pp.~1075--1077.

\bibitem{dimassi1999spectral}
{\sc M.~Dimassi and J.~Sj{\"o}strand}, {\em Spectral asymptotics in the
  semi-classical limit}, no.~268, Cambridge university press, 1999.

\bibitem{dombrowski2011quantization}
{\sc N.~Dombrowski, F.~Germinet, and G.~Raikov}, {\em Quantization of edge
  currents along magnetic barriers and magnetic guides}, in Annales Henri
  Poincar{\'e}, vol.~12, Springer, 2011, pp.~1169--1197.

\bibitem{drouot2022semiclassical}
{\sc A.~Drouot}, {\em Topological insulators in semiclassical regime}, Arxiv
  preprint arXiv:2206.08238.

\bibitem{Drouot:19b}
\leavevmode\vrule height 2pt depth -1.6pt width 23pt, {\em The bulk-edge
  correspondence for continuous honeycomb lattices}, Communication in Partial
  Differential Equations, 44 (2019), p.~1406–1430.

\bibitem{Drouot2020microlocal}
\leavevmode\vrule height 2pt depth -1.6pt width 23pt, {\em Microlocal analysis
  of the bulk-edge correspondence}, Communications in Mathematical Physics, 383
  (2021), p.~2069–2112.

\bibitem{Drouot:21}
\leavevmode\vrule height 2pt depth -1.6pt width 23pt, {\em Ubiquity of conical
  points in topological insulators}, Journal de l'Ecole Polytechnique, 8
  (2021), pp.~507--532.

\bibitem{drouot2020edge}
{\sc A.~Drouot and M.~Weinstein}, {\em Edge states and the valley hall effect},
  Advances in Mathematics, 368 (2020), p.~107142.

\bibitem{EG02}
{\sc P.~Elbau and G.-M. Graf}, {\em Equality of bulk and edge {H}all
  conductance revisited}, Communications in mathematical physics, 229 (2002),
  pp.~415--432.

\bibitem{GP}
{\sc G.-M. {Graf} and M.~{Porta}}, {\em Bulk-edge correspondence for
  two-dimensional topological insulators}, Communications in Mathematical
  Physics, 324 (2013), pp.~851--895.

\bibitem{GS-CUP-94}
{\sc A.~Grigis and j.~Sj{\"o}strand}, {\em Microlocal Analysis for Differential
  Operators: An Introduction}, Cambridge University Press, 1994.

\bibitem{guillemin1990geometric}
{\sc V.~Guillemin and S.~Sternberg}, {\em Geometric asymptotics}, no.~14,
  American Mathematical Soc., 1990.

\bibitem{hall2013quantum}
{\sc B.~C. Hall}, {\em Quantum theory for mathematicians}, vol.~267, Springer,
  2013.

\bibitem{hu2022traveling}
{\sc P.~Hu, P.~Xie, and Y.~Zhu}, {\em Traveling edge states in massive dirac
  equations along slowly varying edges}, arXiv preprint arXiv:2202.13653,
  (2022).

\bibitem{lee2019elliptic}
{\sc J.~P. Lee-Thorp, M.~I. Weinstein, and Y.~Zhu}, {\em Elliptic operators
  with honeycomb symmetry: Dirac points, edge states and applications to
  photonic graphene}, Archive for Rational Mechanics and Analysis, 232 (2019),
  pp.~1--63.

\bibitem{lindblad2020modified}
{\sc H.~Lindblad, J.~Luhrmann, W.~Schlag, and A.~Soffer}, {\em On modified
  scattering for 1d quadratic klein-gordon equations with non-generic
  potentials}, arXiv preprint arXiv:2012.15191,  (2020).

\bibitem{moessner2021topological}
{\sc R.~Moessner and J.~E. Moore}, {\em Topological Phases of Matter},
  Cambridge University Press, 2021.

\bibitem{paul1993construction}
{\sc T.~Paul and A.~Uribe}, {\em A construction of quasi-modes using coherent
  states}, in Annales de l'IHP Physique th{\'e}orique, vol.~59, 1993,
  pp.~357--381.

\bibitem{PSB16}
{\sc E.~Prodan and H.~Schulz-Baldes}, {\em Bulk and boundary invariants for
  complex topological insulators}, Springer verlag, Berlin, 2016.

\bibitem{QB-NUMTI-2021}
{\sc S.~Quinn and G.~Bal}, {\em {Approximations of interface topological
  invariants}}, arXiv:2112.02686,  (2022).

\bibitem{quinn2022asymmetric}
\leavevmode\vrule height 2pt depth -1.6pt width 23pt, {\em Asymmetric transport
  for magnetic dirac equations}, arXiv preprint arXiv:2211.00726,  (2022).

\bibitem{ramond1996semiclassical}
{\sc T.~Ramond}, {\em Semiclassical study of quantum scattering on the line},
  Communications in mathematical physics, 177 (1996), pp.~221--254.

\bibitem{souslov2019topological}
{\sc A.~Souslov, K.~Dasbiswas, M.~Fruchart, S.~Vaikuntanathan, and V.~Vitelli},
  {\em Topological waves in fluids with odd viscosity}, Physical Review
  Letters, 122 (2019), p.~128001.

\bibitem{thaller2013dirac}
{\sc B.~Thaller}, {\em {The Dirac equation}}, Springer Science \& Business
  Media, 2013.

\bibitem{TKNN}
{\sc D.~J. Thouless, M.~Kohmoto, M.~P. Nightingale, and M.~den Nijs}, {\em
  {Quantized Hall Conductance in a Two-Dimensional Periodic Potential}}, Phys.
  Rev. Lett., 49 (1982), pp.~405--408.

\bibitem{volovik2009universe}
{\sc G.~Volovik}, {\em The Universe in a Helium Droplet}, International Series
  of Monographs on Physics, OUP Oxford, 2009.

\bibitem{watanabe2020counting}
{\sc H.~Watanabe}, {\em Counting rules of nambu--goldstone modes}, Annual
  Review of Condensed Matter Physics, 11 (2020), pp.~169--187.

\bibitem{WI}
{\sc E.~{Witten}}, {\em {Three lectures on topological phases of matter}},
  Nuovo Cimento Rivista Serie, 39 (2016), pp.~313--370.

\bibitem{zworski2012semiclassical}
{\sc M.~Zworski}, {\em Semiclassical analysis}, vol.~138, American Mathematical
  Soc., 2012.

\end{thebibliography}

}

%
\end{document}